\documentclass[11pt]{article}

\usepackage[T1]{fontenc}
\usepackage{lmodern}
\usepackage{etoolbox}

\makeatletter
\patchcmd{\@maketitle}{\LARGE \@title}{\LARGE\bfseries\@title}{}{}
\makeatother

\usepackage{amsmath, amssymb, amsthm} 
\usepackage[margin=2.5cm]{geometry}
\usepackage{color, graphicx, soul}
\usepackage{caption}

\usepackage{hyperref}
\hypersetup{
	colorlinks=true,        
	linkcolor=black,         
	citecolor=black,         
	urlcolor=black           
}

\usepackage[overload]{empheq}
\definecolor{myblue}{rgb}{.9, .9, 1}

\usepackage[capitalize]{cleveref}
\crefname{equation}{}{equations}
\crefname{chapter}{Appendix}{chapters}
\crefname{item}{}{items}
\crefname{figure}{Figure}{figures}
\makeatletter
\def\namedlabel#1#2{\begingroup
	\def\@currentlabel{#2}%
	\label{#1}\endgroup
}
\makeatother

\makeatletter
\def\th@plain{%
	\thm@notefont{}
	\itshape 
}
\def\th@definition{%
	\thm@notefont{}
	\normalfont 
}

\renewenvironment{proof}[1][\proofname]{\par
	\normalfont
	\topsep0\p@\@plus3\p@ \trivlist
	\item[\hskip\labelsep\itshape
	#1\@addpunct{.}]\ignorespaces
}{%
	\qed\endtrivlist
}
\makeatother

\newtheorem{theorem}{Theorem}[section]
\newtheorem{lemma}[theorem]{Lemma}
\newtheorem{corollary}[theorem]{Corollary}
\newtheorem{proposition}[theorem]{Proposition}
\newtheorem{fact}[theorem]{Fact}

\theoremstyle{definition}
\newtheorem{definition}[theorem]{Definition}
\theoremstyle{definition}
\newtheorem{example}[theorem]{Example}
\theoremstyle{definition}
\newtheorem{remark}[theorem]{Remark}

\usepackage[shortlabels]{enumitem}
\setlist[enumerate]{nosep}

\allowdisplaybreaks

\newcommand{\menge}[2]{\{{#1}~\big |~{#2}\}} 
\newcommand{\mmenge}[2]{\bigg\{{#1}~\bigg |~{#2}\bigg\}}

\newcommand{\ball}[2]{{\ensuremath{\it I\hspace{-5pt}B}}(#1;#2)}
\newcommand{\bball}[2]{{\ensuremath{\it I\hspace{-5pt}B}}\big(#1;#2\big)}

\newcommand{\scal}[2]{\left\langle {#1},{#2} \right\rangle}

\newcommand{\pnX}[1]{N^{\rm prox}_{#1}} 

\newcommand{\To}{\ensuremath{\rightrightarrows}}

\newcommand{\NN}{\ensuremath{\mathbb N}}
\newcommand{\nnn}{\ensuremath{{n\in{\mathbb N}}}}

\newcommand{\RR}{\ensuremath{\mathbb R}}
\newcommand{\RP}{\ensuremath{\mathbb{R}_+}}
\newcommand{\RPP}{\ensuremath{\mathbb{R}_{++}}}

\newcommand{\argmin}{\ensuremath{\operatorname*{argmin}}}

\newcommand{\reli}{\ensuremath{\operatorname{ri}}}

\newcommand{\aff}{\ensuremath{\operatorname{aff}}}

\newcommand{\Fix}{\ensuremath{\operatorname{Fix}}}
\newcommand{\Id}{\ensuremath{\operatorname{Id}}}

\def\ox{\overline{x}}
\def\oy{\overline{y}}
\def\O{\Omega}

\def\disp{\displaystyle}

\begin{document}

\title{Linear Convergence of Projection Algorithms}

\author{
Minh N.\ Dao\thanks{CARMA, University of Newcastle, Callaghan, NSW 2308, Australia. 
E-mail: \texttt{daonminh@gmail.com}}
~and~
Hung M.\ Phan\thanks{Department of Mathematical Sciences, Kennedy College of Sciences, University of Massachusetts  Lowell, MA 01854, USA.
E-mail: \texttt{hung\char`_phan@uml.edu}.}}

\date{August 16, 2017}

\maketitle

\begin{abstract}
Projection algorithms are well known for their simplicity and flexibility in solving feasibility problems. They are particularly important in practice due to minimal requirements for software implementation and maintenance. In this work, we study linear convergence of several projection algorithms for systems of finitely many closed sets. The results complement contemporary research on the same topic.
\end{abstract}

{\small

}

\section{Introduction}
In this paper, $X$ is a Euclidean space with inner product $\scal{\cdot}{\cdot}$ and induced norm $\|\cdot\|$. Throughout, we set $I :=\{1,\dots,m\}$ and assume that $\{C_i\}_{i\in I}$ is a system of closed (possibly nonconvex) subsets of $X$. The notation used in the paper is fairly standard and follows \cite{BC11}. The nonnegative integers are $\NN$, the real numbers are $\RR$, while $\RP := \menge{x \in \RR}{x \geq 0}$ and $\RPP := \menge{x \in \RR}{x >0}$.
If $w \in X$ and $\rho \in \RP$, then $\ball{w}{\rho} :=\menge{x \in X}{\|x -w\| \leq \rho}$ is the closed ball centered at $w$ with radius $\rho$. Given a subset $C$ of $X$, the affine hull of $C$ is denoted by $\aff C$ and the orthogonal complement of $C$ is $C^\perp :=\menge{x\in X}{\forall c\in C:\ \scal{c}{x} =0}$. The notation $T:X\rightrightarrows X$ means that $T$ is a set-valued operator from $X$ to $X$ and $\Fix T :=\menge{x\in X}{x\in Tx}$ denotes the set of fixed points of $T$. As usual, $\Id$ represents the identity operator. 

The paper is concerned with \emph{cyclic algorithms} for solving the feasibility problem
\begin{equation}\label{e:feas}
\text{find a point}\quad
x\in \bigcap_{i\in I}C_i.
\end{equation}
This problem has long been known for its importance in many applications. To describe cyclic algorithms for \eqref{e:feas}, we first associate each set $C_i$ with an operator $T_i\colon X\To X$ and adopt the following convention
\begin{equation}
\label{e:cvn}
\forall\nnn,\ \forall i \in I:\quad C_{mn+i}:= C_i
\quad\text{and}\quad
T_{mn+i}:= T_i.
\end{equation}
Given a starting point $x_0 \in X$, the cyclic algorithm with respect to the ordered tuple $(T_i)_{i\in I}$ generates sequences $(x_n)_\nnn$ by
\begin{equation}
\label{e:cycseq}
\forall\nnn:\quad x_{n+1}\in T_{n+1} x_n.
\end{equation}
Each such sequence is called a \emph{cyclic sequence} generated by $(T_i)_{i\in I}$. When $m =1$, we drop the subscripts and write $C :=C_1$ and $T :=T_1$. The recurrence \eqref{e:cycseq} then reads as
\begin{equation}
\forall\nnn:\quad x_{n+1}\in Tx_n,
\end{equation}
and we say that the sequence $(x_n)_\nnn$ is generated by $T$. 
The corresponding operators include, but not limited to, projectors and their variants. Recall that for a set $C$, the \emph{distance function} to $C$ is defined by
\begin{equation}
d_C\colon X\to\RR\colon x\mapsto \inf_{c \in C}\|x -c\|,
\end{equation}
and the \emph{projector} onto $C$ is defined by
\begin{equation}
P_C\colon X\To C\colon x\mapsto \argmin_{c \in C}\|x -c\| =\menge{c \in C}{\|x -c\| =d_C(x)}.
\end{equation}
In general, one expects the cyclic sequence $(x_n)_\nnn$ or other acquired sequences converge to a solution of \eqref{e:feas}. In such case, we are interested in $R$-linear convergence of those sequences. Recall that a sequence $(x_n)_\nnn$ is said to \emph{converge $R$-linearly} to a point $\ox$ with rate $\rho\in\left[0,1\right[$ if there exists a constant $\sigma\in\RP$ such that
\begin{equation}
\forall\nnn:\quad\|x_n-\ox\|\leq\sigma\rho^n.
\end{equation}

Among the main contributions of the paper, under certain regularity assumptions on sets and system of sets,
we show that:
\begin{enumerate}[label = (R\arabic*),leftmargin=0.4in]
\item\label{r1} The cyclic relaxed projections with at most one reflection, which includes the reflection-projection algorithm \cite{BKr04}, converge $R$-linearly locally (see Theorem~\ref{t:linCP} and Remark~\ref{r:1stRP});
\item\label{r2} A {\em refined} $R$-linear rate is obtained for cyclic {over-relaxed} projections (see Theorem~\ref{t:linCP_m-1} and Corollary~\ref{c:cycp});
\item\label{r3} The cyclic semi-intrepid projections for injectable sets converge locally with $R$-linear rate (see Theorem~\ref{t:cycsIP});
\end{enumerate}
Moreover, the linear convergence is global in the presence of convexity (see Corollaries~\ref{c:cvxCP} and \ref{c:cycsIP}). To the best of our knowledge, these results are {\em new} and have not been observed in the literature. In addition, we also present other new results involving Douglas--Rachford (DR) operators \cite{DR56,LM79}; see Theorems~\ref{t:cDR} and \ref{t:0830a}. Our work complements other studies on projection algorithms \cite{BD17,BDNP16a,BDNP16b,BNP15,BPW13,DIL15,HL13,LLM09,NR16,Pha14}.

The remainder of the paper is organized as follows. Section~\ref{s:pre} contains basic concepts needed for our analysis. Section~\ref{s:qFm} then provides key components for $R$-linear convergence. In Section~\ref{s:lincyc}, we prove $R$-linear convergence for general cyclic algorithms. Finally, Section~\ref{s:apps} presents applications to various cyclic algorithms including the cyclic relaxed projections, cyclic semi-intrepid projections, and cyclic generalized DR algorithm.

\section{Preliminaries}
\label{s:pre}
Given a subset $C$ of $X$ and $x\in C$, the \emph{Fr\'echet normal cone} to $C$ at $x$ \cite[Definition~1.1(i)]{Mor06} is defined by 
\begin{equation}
\widehat{N}_C(x) :=\mmenge{u\in X}{\limsup_{y\to x,\, y\in C\smallsetminus\{x\}} \frac{\scal{u}{y -x}}{\|y -x\|} \leq 0},
\end{equation}
the \emph{proximal normal cone} to $C$ at $x$ (see \cite[Section~2.5.2, D]{Mor06} and \cite[Example~6.16]{RW98}) is given by
\begin{equation}
\pnX{C}(x) :=\menge{\lambda(z-x)}{z\in P_C^{-1}(x),\ \lambda\in\RP},
\end{equation}
and the \emph{limiting normal cone} to $C$ at $x$ \cite[Definition~1.1(ii)]{Mor06} can be given by (\cite[Theorem~1.6]{Mor06})
\begin{subequations}\label{e:170624b}
\begin{align}
N_C(x) &:=\menge{u \in X}{\exists x_n\to x, u_n\to u \text{~with~} x_n\in C, u_n\in \widehat{N}_C(x_n)}\\
&\phantom{:}=\menge{u \in X}{\exists x_n\to x, u_n\to u \text{~with~} x_n\in C, u_n\in \pnX{C}(x_n)}.
\end{align}
\end{subequations}
As seen below, normal cones are used to describe superregularity for sets and strong regularity for systems of sets. We recall the superregularity concept, which was first introduced in \cite{LLM09} and later refined in \cite{BLPW13a,BLPW13b,HL13,NR16}. Superregularity holds for a major class of sets including convex sets and sets with ``smooth" boundary. This property plays an important role in analyzing linear convergence of projection methods, see, e.g., \cite{BLPW13a,BLPW13b,HL13,LLM09,NR16,Pha14}.

\begin{definition}[superregularity of sets]
\label{d:supreg}
Let $C$ be a nonempty subset of $X$, $w\in X$, $\varepsilon \in \RP$, and $\delta \in \RPP$.
We say that $C$ is \emph{$(\varepsilon,\delta)$-regular} at $w$ if
\begin{equation}
\left.\begin{aligned}
&x,y\in C\cap\ball{w}{\delta},\\
&u\in\pnX{C}(x)
\end{aligned}\right\}\ \Rightarrow\
\scal{u}{x-y}\geq -\varepsilon\|u\|\cdot\|x-y\|,
\end{equation}
and \emph{$(\varepsilon, \infty)$-regular} at $w$ if it is $(\varepsilon, \delta)$-regular for all $\delta\in \RPP$. The set $C$ is said to be \emph{superregular} at $w$ if for all $\varepsilon \in \RPP$, there exists $\delta \in \RPP$ such that $C$ is $(\varepsilon,\delta)$-regular at $w$. The system $\{C_i\}_{i\in I}$ is said to be \emph{superregular} at $w$ if $C_i$ is superregular at $w$ for every $i\in I$.
\end{definition}

Next, we recall two regularity concepts for systems of sets: linear regularity and strong regularity.

\begin{definition}[linear regularity of set systems]
\label{d:linreg}
Let $\kappa \in\RPP$. The system $\{C_i\}_{\in I}$ is said to be \emph{$\kappa$-linearly regular} on a subset $U$ of $X$ if 
\begin{equation}
\label{e:linreg}
\forall x\in U:\quad d_C(x) \leq \kappa\max_{i\in I}d_{C_i}(x),
\quad\text{where}\quad C :=\bigcap_{i\in I}C_i.
\end{equation}
The constant $\kappa$ is called a \emph{linear regularity modulus} of $\{C_i\}_{i\in I}$ on $U$.
We say that $\{C_i\}_{i\in I}$ is \emph{linearly regular} around $w \in X$ if there exist $\delta \in \RPP$ and $\kappa\in \RPP$ such that $\{C_i\}_{i\in I}$ is $\kappa$-linearly regular on $\ball{w}{\delta}$.
\end{definition}

Linear regularity for set systems has a long history and was first defined in convex settings, see, e.g., \cite[Definition~5.6]{BB96}, \cite[Definition~3]{BBL99},
and \cite[Section~5.2]{BNP15} for a brief summary on this property. Naturally, linear regularity was extended to system of closed sets, for instance, \cite[Definition~3.5]{HL13}; and was known as \emph{metric inequality} in \cite[Equation~(15)]{Kru06}, \cite[Section~3]{NT01}, and \cite[Section~5]{Iof89}; and as {\em subtransversality} in \cite[Definition~1]{KLN17}.

\begin{definition}[strong regularity of set systems]
\label{d:strgreg}
The system $\{C_i\}_{\in I}$ is said to be \emph{strongly regular} at $w\in \bigcap_{i\in I} C_i$ if 
\begin{equation}\label{e:strgreg}
\sum_{i\in I} u_i =0 \text{ \ and \ } u_i\in N_{C_i}(w)
\quad\Rightarrow\quad \forall i\in I:\ u_i=0.
\end{equation}
In the case $I =\{1,2\}$, condition \eqref{e:strgreg} can be rewritten as
\begin{equation}
N_{C_1}(w)\cap (-N_{C_2}(w)) =\{0\}.
\end{equation}
\end{definition}

Strong regularity of systems is also known as \emph{normal qualification condition} in \cite[Definition~3.2]{Mor06}, as \emph{CQ condition} in \cite[Definition~6.6]{BLPW13a}, and as {\em transversality} in \cite[Definition~2]{KLN17}. To clear the confusion it may cause, we will show that strong regularity in Definition~\ref{d:strgreg} is equivalent to the ones in \cite[Definition~1(vi)]{Kru06} and in \cite[Definition~3.2]{HL13}. In view of \cite[Proposition~2, Proposition~10(ii), and Corollary~2]{Kru06}, it suffices to prove the following result.

\begin{proposition}[characterization of strong regularity]
The system $\{C_i\}_{i\in I}$ is strongly regular at $w\in \bigcap_{i\in I}C_i$ if and only if there exist $\zeta\in \RPP$ and $\delta\in\RPP$ such that
\begin{equation}\label{e:170629b}
\forall i\in I, \forall x_i\in C_i\cap\ball{w}{\delta}, \forall u_i\in \widehat{N}_{C_i}(x_i):\quad 
\Big\|\sum_{i\in I}u_i\Big\|\geq \zeta\sum_{i\in I}\|u_i\|.
\end{equation}
\end{proposition}
\proof 
$(\Leftarrow)$: Suppose that \eqref{e:170629b} holds and that $\sum_{i\in I} u_i =0$ with $u_i\in N_{C_i}(w)$. Then for every $i\in I$, by \eqref{e:170624b}, there exist sequences $x_{i,n}\to w$, $u_{i,n}\to u_i$ with $x_{i,n}\in C_i$ and $u_{i,n}\in \widehat{N}_{C_i}(x_{i,n})$. 
Since $x_{i,n}\to w$, we can assume without loss of generality that $x_{i,n}\in \ball{w}{\delta}$ for all $\nnn$. It follows that $x_{i,n}\in C_i\cap \ball{w}{\delta}$, and then by  \eqref{e:170629b}, we have
$\big\|\sum_{i\in I}u_{i,n}\big\|\geq \zeta\sum_{i\in I}\|u_{i,n}\|$ for all $\nnn$.
Passing to the limit as $n\to\infty$, we get $\|\sum_{i\in I}u_i\|\geq \zeta\sum_{i\in I}\|u_i\|$. Combining with the assumption $\sum_{i\in I} u_i =0$, we derive $u_i =0$ for every $i\in I$. 
 
$(\Rightarrow)$: Suppose to the contrary that \eqref{e:170629b} is not true. Then there exist sequences $\zeta_n\to 0^+$, $\delta_n\to 0^+$, $x_{i,n}\in C_i\cap\ball{w}{\delta_n}$, and $u_{i,n}\in \widehat{N}_{C_i}(x_{i,n})$ such that 
\begin{equation}\label{e:170629c}
\forall\nnn:\quad \Big\|\sum_{i\in I}u_{i,n}\Big\|<\zeta_n\sum_{i\in I}\|u_{i,n}\|
\quad\text{and}\quad \sum_{i\in I}\|u_{i,n}\| =1,
\end{equation}
where the latter is obtained by rescaling if necessary.
Thus, for every $i\in I$, the sequence $(u_{i,n})_\nnn$ is bounded, and by extracting subsequences, we can assume that $u_{i,n}\to u_i$. Since $x_{i,n}\to w$ and $x_{i,n}\in C_i$,  
it follows from \eqref{e:170624b} that $u_i\in N_{C_i}(w)$. Letting $n\to\infty$ in \eqref{e:170629c}, we obtain
$\big\|\sum_{i\in I}u_{i}\big\|=0$
and $\sum_{i\in I}\|u_i\|=1$,
which contradicts the strong regularity. Thus \eqref{e:170629b} holds.
\endproof

We end this section with connections between linear regularity and strong regularity.
\begin{fact}\label{f:str-lin}
{\rm(\cite[Theorem~1]{Kru06})}
If the system $\{C_i\}_{i\in I}$ is strongly regular at $w \in \bigcap_{i\in I}C_i$, then it is linearly regular around $w$.
\end{fact}

\begin{remark}[strong regularity of subsystems]
\label{r:subsystem}
By definition, if the system $\{C_i\}_{i\in I}$ is strongly regular at $w$, then so is each of its subsystems. 
However, even when each proper subsystem $\{C_i\}_{i\in J}$ with $J\subsetneqq I$ is strongly regular and the entire system $\{C_i\}_{i\in I}$ is linearly regular, it does not imply that $\{C_i\}_{i\in I}$ is strongly regular. 
For example, in $\RR^2$, consider $C_1 =\menge{(\xi, \zeta)}{\xi +\zeta\leq 0}$, $C_2 =\menge{(\xi, \zeta)}{\xi -\zeta\leq 0}$, $C_3 =\menge{(\xi, \zeta)}{\xi\geq 0}$, and $w =(0,0)\in C_1\cap C_2\cap C_3$. Then one can check that $\{C_i\}_{i\in J}$ with $J\subsetneqq\{1,2,3\}$ is strongly regular at $w$, and $\{C_1,C_2,C_3\}$ is linearly regular around $w$, but $\{C_1,C_2,C_3\}$ is not strongly regular at $w$. 
\end{remark}

\section{Quasi Fej\'er monotonicity and quasi coercivity}
\label{s:qFm}
The following quasi Fej\'er monotonicity concept generalizes the Fej\'er monotonicity for sequences and operators, see, e.g., \cite[Definition~5.1]{BC11} and \cite[Definition~2.1.15]{Ceg12}.

\begin{definition}[quasi firm Fej\'er monotonicity]
\label{d:qfF}
Let $C$ and $U$ be nonempty subsets of $X$, let $\gamma \in \left[1, +\infty\right[$, and let $\beta \in \RP$. A set-valued operator $T\colon X\To X$ is said to be \emph{$(C, \gamma, \beta)$-quasi firmly Fej\'er monotone} on $U$ if
\begin{equation}
\label{e:qfF}
\forall x \in U,\ \forall x_+ \in Tx,\ \forall \ox \in C:\quad
\|x_+ -\ox\|^2 +\beta\|x -x_+\|^2 \leq \gamma \|x -\ox\|^2.
\end{equation}
We say that $T$ is \emph{$(C, \gamma)$-quasi firmly Fej\'er monotone} on $U$ if $\beta=1$, i.e.,
\begin{equation}
\forall x \in U,\ \forall x_+ \in Tx,\ \forall \ox \in C:\quad
\|x_+ -\ox\|^2 +\|x -x_+\|^2 \leq \gamma \|x -\ox\|^2,
\end{equation}
and that $T$ is \emph{$(C, \gamma)$-quasi Fej\'er monotone} on $U$ if $\beta=0$, i.e.,
\begin{equation}
\forall x \in U,\ \forall x_+ \in Tx,\ \forall \ox \in C:\quad
\|x_+ -\ox\| \leq \gamma^{1/2} \|x -\ox\|.
\end{equation}
\end{definition}

From the definition, we observe that
\begin{enumerate}
\item $(C, \gamma, \beta)$-quasi firm Fej\'er monotonicity implies $(C, \gamma)$-quasi Fej\'er monotonicity, while $(C, 1)$-quasi Fej\'er monotonicity is exactly Fej\'er monotonicity with respect to $C$ in \cite[Definition~2.1.15]{Ceg12}.	
\item If $\gamma' \geq \gamma \geq 1$, $0 \leq \beta' \leq \beta$, $C' \subseteq C$, and $U' \subseteq U$, then $(C, \gamma, \beta)$-quasi firm Fej\'er monotonicity on $U$ implies $(C', \gamma', \beta')$-quasi firm Fej\'er monotonicity on $U'$.
\item If $T$ is nonexpansive (see \cite[Definition~4.1]{BC11}), then $T$ is $(\Fix T, 1)$-quasi Fej\'er monotone on $X$.
\item If $T$ is $\lambda$-averaged (see \cite[Definition~4.23]{BC11}), then by \cite[Proposition~4.25(iii)]{BC11}, $T$ is $(\Fix T, 1, \frac{1 -\lambda}{\lambda})$-quasi firmly Fej\'er monotone on $X$. In particular, if $T$ is firmly nonexpansive, then $T$ is $(\Fix T, 1)$-quasi firmly Fej\'er monotone on $X$.
\end{enumerate}
Quasi firm Fej\'er monotonicity is closely related to \cite[Definition~2.3]{HL13} and \cite[Proposition~2.4(iii)]{LTT16}. 
Also, $(C,\gamma)$-quasi firm Fej\'er monotonicity is more restrictive than \cite[Definition~2.7]{Pha14} since the latter requires only
\begin{equation}
\forall x \in U,\ \forall x_+ \in Tx,\
\forall \ox \in P_C x:\quad
\|x_+ -\ox\|^2 +\|x -x_+\|^2 \leq \gamma \|x -\ox\|^2.
\end{equation}
Nevertheless, it turns out that quasi firm Fej\'er monotonicity still holds for a broad class of operators, e.g., relaxed projectors for superregular sets (see Proposition~\ref{p:rproj}) and generalized Douglas--Rachford operators for systems of two superregular sets (see Proposition~\ref{p:qfF-DR}).

The next lemma shows the quasi firm Fej\'er monotonicity for averaged-type operators.

\begin{lemma}[averaged quasi firmly Fej\'er monotone operators]
\label{l:averaged}
Let $C$ and $U$ be nonempty subsets of $X$, $\gamma \in \left[1, +\infty\right[$, $\beta\in\RP$, $\lambda \in \left]0, 1+\beta\right]$, and let $S\colon X\To X$ be a $(C,\gamma,\beta)$-quasi firmly Fej\'er monotone operator on $U$. Then $T :=(1 -\lambda)\Id +\lambda S$ is $(C,\gamma',\beta')$-quasi firmly Fej\'er monotone on $U$ with
\begin{equation}
\gamma':=1 -\lambda +\lambda\gamma
\quad\text{and}\quad
\beta':=\frac{1-\lambda+\beta}{\lambda}.
\end{equation}
\end{lemma}
\proof 
Let $x \in U$, $x_+ \in Tx$, and $\ox \in C$. Writing $x_+ =(1 -\lambda)x +\lambda s$ with $s \in Sx$, we have $x_+ -\ox =(1 -\lambda)(x -\ox) +\lambda(s -\ox)$ and $x -x_+ =\lambda(x -s)$. So
\begin{subequations}\label{e:160722a}
\begin{align}
\|x_+ -\ox\|^2 &=(1 -\lambda)\|x -\ox\|^2 +\lambda\|s -\ox\|^2 -\lambda(1 -\lambda)\|(x -\ox) -(s -\ox)\|^2 \\
&=(1 -\lambda)\|x -\ox\|^2 +\lambda\|s -\ox\|^2 -\lambda(1 -\lambda)\|x -s\|^2.
\end{align}
\end{subequations} 
Using the $(C,\gamma,\beta)$-quasi firm Fej\'er monotonicity of $S$ on $U$, we continue \eqref{e:160722a} as
\begin{subequations}
\begin{align}
\|x_+ -\ox\|^2
&\leq (1 -\lambda)\|x -\ox\|^2 
+\lambda\Big(\gamma\|x -\ox\|^2-\beta\|x-s\|^2\Big)-\lambda(1 -\lambda)\|x -s\|^2\\ 
&=(1 -\lambda +\lambda\gamma)\|x -\ox\|^2-\lambda(1-\lambda+\beta)\|x-s\|^2\\
&=(1 -\lambda +\lambda\gamma)\|x -\ox\|^2-\frac{1-\lambda+\beta}{\lambda}\|x-x_+\|^2.
\end{align}
\end{subequations}
This completes the proof.
\endproof

\begin{definition}[quasi coercivity]
\label{d:qco}
Let $C$ and $U$ be nonempty subsets of $X$ and let $\nu \in \RPP$. An operator $T\colon X\To X$ is said to be \emph{$(C, \nu)$-quasi coercive} on $U$ if
\begin{equation}
\forall x\in U,\ \forall x_+\in Tx:\quad \|x -x_+\|\geq \nu d_C(x).
\end{equation}
We say that $T$ is \emph{$C$-quasi coercive} around $w \in X$ if there exist $\delta\in\RPP$ and $\nu\in \RPP$ such that $T$ is $(C, \nu)$-quasi coercive on $\ball{w}{\delta}$.
\end{definition}
Obviously, if $0< \nu'\leq \nu$, $C'\supseteq C$, and $U'\subseteq U$, then $(C, \nu)$-quasi coercivity on $U$ implies $(C', \nu')$-quasi coercivity on $U'$. Quasi coercivity follows and slightly extends the {\em coercivity condition} in \cite[Lemma~3.1(b)]{HL13} because the latter requires $C\subseteq \Fix T$ while the former does not.
Quasi coercivity is also closely related to the {\em linear regularity for operators} in \cite[Definition~2.1]{BNP15}. Indeed, when $C =\Fix T \neq\varnothing$, then $T$ is $(C, \nu)$-quasi coercive on $X$ if and only if it is linearly regular with constant $\frac{1}{\nu}$ in the sense of \cite[Definition~2.1]{BNP15}. Again, under certain conditions, we will show that quasi coercivity holds for several class of projectors.

\subsection{Relaxed projectors}
In this section, we show the quasi firm Fej\'er monotonicity and quasi coercivity of relaxed projectors for superregular sets. Let $C$ be a nonempty closed subset of $X$ and let $\lambda\in \RP$. The \emph{relaxed projector} for $C$ with parameter $\lambda$ is defined by
\begin{equation}
P_C^\lambda :=(1- \lambda)\Id +\lambda P_C.
\end{equation}
We say that $P_C^\lambda$ is \emph{under-relaxed} if $\lambda\leq 1$ and \emph{over-relaxed} if $\lambda\geq 1$.
Clearly, $P_C^0 =\Id$, $P_C^1 =P_C$, and $P_C^2 =R_C :=2P_C -\Id$ (the \emph{reflector} across $C$). 
The following lemma will be used several times in our analysis.

\begin{lemma}
\label{l:ImBall}
Let $w \in C$, let $\gamma \in \left[1, +\infty\right[$, and let $\delta \in \RPP$. Then the following hold:
\begin{enumerate}
\item\label{l:ImBall_P}
For all $x\in\ball{w}{\delta/2}$, $P_C^\lambda x \subseteq \ball{w}{(1 +\lambda)\delta/2}$. In particular, 
$P_C(\ball{w}{\delta/2}) \subseteq C\cap \ball{w}{\delta}$.
\item\label{l:ImBall_T}
If $T\colon X\To X$ is $(C\cap \ball{w}{\delta}, \gamma)$-quasi Fej\'er monotone on $\ball{w}{\delta/2}$, then
\begin{subequations}
\begin{align}
\forall x \in \ball{w}{\delta/2}:&\quad Tx \subseteq \ball{w}{\gamma^{1/2}\delta/{2}}, \\
\forall x \in \ball{w}{\delta/2},\ \forall x_+ \in Tx:&\quad d_C(x_+) \leq \gamma^{1/2}d_C(x).
\end{align}
\end{subequations} 
\end{enumerate}
\end{lemma}
\proof
\ref{l:ImBall_P}: Let $x \in \ball{w}{\delta/2}$ and let $x_+ \in P_C^\lambda x$. Writing $x_+ =(1 -\lambda)x +\lambda p$ for some $p \in P_C x$ and noting that $w \in C$, we have $\|x_+ -x\| =\lambda\|p -x\| =\lambda d_C(x) \leq \lambda\|x -w\|$ and so
\begin{equation}
\|x_+ -w\| \leq \|x_+ -x\| +\|x -w\| \leq (1 +\lambda)\|x -w\| \leq (1 +\lambda)\delta/2.
\end{equation}
Therefore, $P_C^\lambda x \subseteq C\cap 
\ball{w}{(1+\lambda)\delta/2}$.

\ref{l:ImBall_T}: Let $x \in \ball{w}{\delta/2}$ and let $x_+ \in Tx$. By quasi Fej\'er monotonicity, 
\begin{equation}
\label{e:qF}
\forall\ox \in C\cap \ball{w}{\delta}:\quad \|x_+ -\ox\| \leq \gamma^{1/2}\|x -\ox\|.
\end{equation}
Setting $\ox =w$, we have
$\|x_+ -w\| \leq \gamma^{1/2}\|x -w\| \leq \gamma^{1/2}\delta/2$. 
Hence, $Tx \subseteq \ball{w}{\gamma^{1/2}\delta/2}$.
Now let $p \in P_Cx$. Then $p\in C\cap \ball{w}{\delta}$ by \ref{l:ImBall_P}. Applying \eqref{e:qF} to $\ox =p$ yields
\begin{equation}
d_C(x_+) \leq \|x_+ -p\| \leq \gamma^{1/2}\|x -p\| =\gamma^{1/2}d_C(x).
\end{equation}
\endproof

\begin{proposition}[quasi firm Fej\'er monotonicity of relaxed projectors]
\label{p:rproj}
Let $w \in C$, $\varepsilon \in \left[0, 1\right[$, $\delta \in \RPP$, and $\lambda\in \left]0, 2\right]$. Set 
\begin{equation}
\label{e:setting}
\Omega :=C\cap \ball{w}{\delta},\quad
\gamma :=1 +\frac{\lambda\varepsilon}{1 -\varepsilon},
\quad\text{and}\quad \beta :=\frac{2 -\lambda}{\lambda}.
\end{equation}
Suppose that $C$ is $(\varepsilon, \delta)$-regular at $w$. Then $P_C^\lambda$ is $(\Omega, \gamma, \beta)$-quasi firmly Fej\'er monotone and, in particular, $R_C$ is $(\Omega, \frac{1+\varepsilon}{1 -\varepsilon})$-quasi Fej\'er monotone on $\ball{w}{\delta/2}$.
Additionally, if $\varepsilon \in \left[0, 1/3\right]$, then
\begin{equation}
\forall x \in \ball{w}{\delta/2}:\quad P_C^\lambda x \subseteq \ball{w}{\delta/\sqrt{2}}.
\end{equation}	
\end{proposition}
\proof 
Let $x \in \ball{w}{\delta/2}$ and let $p\in P_C x$. Then $p\in \Omega$ by Lemma~\ref{l:ImBall}\ref{l:ImBall_P}. Since $C$ is $(\varepsilon,\delta)$-regular at $w$ and $x -p \in \pnX{C}(p)$, we have
\begin{equation}
\forall \ox\in\O:\quad
\scal{x -p}{p -\ox}\geq -\varepsilon \|x -p\|\cdot\|p -\ox\|
\geq -\frac{\varepsilon}{2}\left(\|x -p\|^2 +\|p -\ox\|^2\right).
\end{equation}
It then follows that
\begin{subequations}
\begin{align}
\forall \ox\in\O:\quad
\|x-\ox\|^2&=\|x-p\|^2+\|p-\ox\|^2+2\scal{x-p}{p-\ox}\\
&\geq \|x-p\|^2+\|p-\ox\|^2-\varepsilon\left(\|x -p\|^2 +\|p -\ox\|^2\right)\\
&=(1-\varepsilon)\big(\|x-p\|^2+\|p-\ox\|^2\big).
\end{align}
\end{subequations}
So
\begin{equation}
\forall \ox\in\O:\quad
\frac{1}{1-\varepsilon}\|x-\ox\|^2\geq \|x -p\|^2 +\|p -\ox\|^2,
\end{equation}
i.e., $P_C$ is $(\Omega, \frac{1}{1-\varepsilon}, 1)$-quasi firmly Fej\'er monotone on $\ball{w}{\delta/2}$. Now by Lemma~\ref{l:averaged}, we conclude that $P^\lambda_C=(1-\lambda)\Id+\lambda P_C$ is $(\Omega, \gamma, \beta)$-quasi firmly Fej\'er monotone on $\ball{w}{\delta/2}$ with $\gamma$ and $\beta$ given by \eqref{e:setting}.
For $\lambda =2$, we have that $\gamma =\frac{1+\varepsilon}{1-\varepsilon}$, $\beta =0$, and so $R_C =P_C^2$ is $(\Omega, \frac{1+\varepsilon}{1 -\varepsilon})$-quasi Fej\'er monotone on $\ball{w}{\delta/2}$.

Next assume that $\varepsilon \in \left[0, 1/3\right]$. Then $\gamma =1 -\lambda +\frac{\lambda}{1 -\varepsilon} \leq 1+\frac{\lambda}{2} \leq 2$. By quasi Fej\'er monotonicity and Lemma~\ref{l:ImBall}\ref{l:ImBall_T}, 
$P_C^\lambda x \subseteq \ball{w}{\gamma^{1/2}\delta/2} \subseteq \ball{w}{\delta/\sqrt{2}}$.
\endproof

\begin{proposition}[quasi coercivity of relaxed projectors]
\label{p:coer-rproj}
If $\lambda \in \RPP$, then $P_C^\lambda$ is $(C, \lambda)$-quasi coercive on $X$.
\end{proposition}
\proof 
Let $x \in X$ and let $x_+ \in P_C^\lambda x$. 
Then $x_+ =(1 -\lambda)x +\lambda p$ for some $p \in P_Cx$. So
$\|x -x_+\| =\lambda\|x -p\| =\lambda d_C(x)$.
\endproof

\subsection{Generalized Douglas--Rachford operators}

In this section, we establish the quasi firm Fej\'er monotonicity and quasi coercivity of generalized Douglas--Rachford operators for systems of two superregular sets.
Let $A$ and $B$ be closed subsets of $X$ such that $A\cap B \neq\varnothing$ and let $\lambda, \mu, \alpha\in \RPP$. 
The {\em generalized Douglas--Rachford operator} for $(A, B)$ with parameters $(\lambda,\mu,\alpha)$ is defined by
\begin{equation}\label{e:170729a}
T_{\lambda, \mu}^\alpha :=(1 -\alpha)\Id +\alpha P_B^\mu P_A^\lambda.
\end{equation}
Note that $T_{1,1}^1=P_BP_A$ is the classical alternating projection operator \cite{Bre65} and that $T_{2,2}^{1/2}=\frac{1}{2}(\Id +R_BR_A)$ is the classical DR operator \cite{DR56,LM79}.

\begin{proposition}[quasi firm Fej\'er monotonicity of generalized DR operators]
\label{p:qfF-DR}
Let $w \in A\cap B$, $\varepsilon_1 \in \left[0, 1/3\right]$, $\varepsilon_2 \in \left[0, 1\right[$, $\delta \in \RPP$, $\lambda, \mu \in \left]0, 2\right]$, and $\alpha \in \left]0, 1\right]$. 
Suppose that $A$ and $B$ are $(\varepsilon_1, \delta)$- and $(\varepsilon_2, \sqrt{2}\delta)$-regular at $w$, respectively.
Then $T_{\lambda, \mu}^\alpha$ is $(A\cap B\cap \ball{w}{\delta}, \gamma, \beta)$-quasi firmly Fej\'er monotone on $\ball{w}{\delta/2}$ with
\begin{equation}
\gamma :=1 -\alpha +\alpha\left(1 +\frac{\lambda\varepsilon_1}{1 -\varepsilon_1}\right)\left(1 +\frac{\mu\varepsilon_2}{1 -\varepsilon_2}\right)
\quad\text{and}\quad \beta :=\frac{1-\alpha}{\alpha}.
\end{equation}
\end{proposition}
\proof 
Let $x \in \ball{w}{\delta/2}$, let $r \in P_A^\lambda x$, let $s \in P_B^\mu r$, and let $\ox \in A\cap B\cap \ball{w}{\delta}$.
Then Proposition~\ref{p:rproj} applied to $P_A^\lambda$ yields
\begin{equation}
\label{e:rPA}
\|r -\ox\| \leq \gamma_1^{1/2}\|x -\ox\|, \quad\text{where}\quad \gamma_1 :=1 +\frac{\lambda\varepsilon_1}{1 -\varepsilon_1},
\end{equation}
and also $r \in \ball{w}{\delta/\sqrt{2}}$.
Next, Proposition~\ref{p:rproj} applied to $P_B^\mu$ yields
\begin{equation}
\|s -\ox\| \leq \gamma_2^{1/2}\|r -\ox\| \leq (\gamma_1\gamma_2)^{1/2}\|x -\ox\|,
\quad\text{where}\quad \gamma_2 :=1 +\frac{\mu\varepsilon_2}{1 -\varepsilon_2}.
\end{equation}
This proves $(A\cap B\cap \ball{w}{\delta}, \gamma_1\gamma_2)$-quasi Fej\'er monotonicity of $P_B^\mu P_A^\lambda$ on $\ball{w}{\delta/2}$. Now apply Lemma~\ref{l:averaged} to the operators $P_B^\mu P_A^\lambda$ and $T_{\lambda,\mu}^\alpha = (1-\alpha)\Id+\alpha P_B^\mu P_A^\lambda$.
\endproof

\begin{proposition}[quasi coercivity of generalized DR operators]
\label{p:coer-DR}
Let $w \in A\cap B$,  
$\lambda,\mu\in \left]0, 2\right]$, and $\alpha \in \RPP$. 
Suppose that $A$ is superregular at $w$ and that $\{A, B\}$ is strongly regular at $w$. 
Then 
\begin{equation}
\overline\theta :=\sup\menge{\scal{u}{v}}{u \in N_A(w)\cap \ball{0}{1},\ v \in \big(-N_B(w)\big)\cap \ball{0}{1}}< 1
\end{equation}
and for all $\theta\in \left]\overline\theta, 1\right[$, there exist $\delta\in \RPP$ and $\kappa\in \RPP$ such that $T_{\lambda, \mu}^\alpha$ is $(A\cap B, \nu)$-quasi coercive on $\ball{w}{\delta/2}$ with
\begin{equation}
\nu:=\frac{\alpha\sqrt{1 -\theta}}{\kappa}\min\Big\{\lambda, \frac{\mu}{\sqrt{1 +\mu^2}}\Big\}.
\end{equation}
\end{proposition}
\proof 
Since $\{A, B\}$ is strongly regular at $w$, we have from \cite[Lemma~2.3]{Pha14} that $\overline\theta< 1$.
Now let $\theta \in \left]\overline\theta, 1\right[$ and let $\varepsilon \in \left[0, 1/3\right]$. 
Using Definition~\ref{d:supreg}, Fact~\ref{f:str-lin}, and \cite[Lemma~4.1]{Pha14}, we can find $\delta\in \RPP$ and $\kappa\in \RPP$ such that $A$ is $(\varepsilon, \delta)$-regular at $w$, that 
\begin{equation}
\label{e:linregAB}
\forall x \in \ball{w}{\delta/2}:\quad
d_{A\cap B}(x) \leq \kappa\max\{d_A(x), d_B(x)\},
\end{equation}
and that 
\begin{equation}
\label{e:strreg}
\left.\begin{aligned}
&a \in A\cap \ball{w}{\delta},\ b \in B\cap \ball{w}{\sqrt{2}\delta}, \\
&u \in \pnX{A}(a),\ v \in \pnX{B}(b)
\end{aligned}\right\}
\quad\Rightarrow\quad
\scal{u}{v}\geq -\theta\|u\|\|v\|.
\end{equation}

Let $x \in \ball{w}{\delta/2}$ and $x_+ \in T_{\lambda, \mu}^\alpha x$. By definition, there exist
\begin{subequations}
\begin{align}[left = \empheqlbrace]
&a \in P_Ax,\quad
r =(1 -\lambda)x +\lambda a\in P_A^\lambda x,\\
&b \in P_Br,\quad
s =(1 -\mu)r +\mu b \in P_B^\mu r
\end{align}
\end{subequations}
such that $x_+ =(1 -\alpha)x +\alpha s$. Then 
\begin{equation}
x -x_+ =\alpha(x -s),\quad x -r =\lambda(x -a), \quad\text{and}\quad r -s =\mu(r -b).
\end{equation}  
By Lemma~\ref{l:ImBall}\ref{l:ImBall_P}, $a \in P_Ax \subseteq A\cap \ball{w}{\delta}$.
Since $\varepsilon \in \left[0, 1/3\right]$, Proposition~\ref{p:rproj} yields $r \in \ball{w}{\sqrt{2}\delta/2}$.
Using again Lemma~\ref{l:ImBall}\ref{l:ImBall_P}, we get $b \in P_Br \subseteq B\cap \ball{w}{\sqrt{2}\delta}$.
Now since $x -r =\lambda(x -a) \in \pnX{A}(a)$ and $r -s =\mu(r -b) \in \pnX{B}(b)$, we use \eqref{e:strreg} and the arithmetic mean-geometric mean inequality to obtain
\begin{equation}
2\scal{x -r}{r -s} \geq -2\theta\|x -r\|\|r -s\| \geq -\theta(\|x -r\|^2 +\|r -s\|^2).
\end{equation}
So
\begin{subequations}\label{e:160802a}
\begin{align}
\|x -x_+\|^2 &=\alpha^2\|x -s\|^2 =\alpha^2(\|x -r\|^2 +\|r -s\|^2 +2\scal{x -r}{r -s}) \\ 
&\geq (1 -\theta)\alpha^2(\|x -r\|^2 +\|r -s\|^2).  
\end{align}
\end{subequations}
Furthermore, 
\begin{equation}\label{e:160802b}
\|x -r\|^2 =\lambda^2\|x -a\|^2 =\lambda^2d_A^2(x),
\end{equation}
and by {the coordinate version of} Cauchy--Schwarz inequality,
\begin{subequations}\label{e:160802c}
\begin{align}
(\mu^2 +1)(\|x -r\|^2 +\|r -s\|^2) &\geq (\mu\|x -r\| +\|r -s\|)^2 \\
&=(\mu\|x -r\| +\mu\|r -b\|)^2 \\
&\geq (\mu\|x -b\|)^2 \geq \mu^2d_B^2(x).
\end{align}
\end{subequations}
Combining \eqref{e:linregAB}, \eqref{e:160802a}, \eqref{e:160802b}, and \eqref{e:160802c}, we obtain
\begin{subequations}
\begin{align}
\|x -x_+\| &\geq \alpha\sqrt{1 -\theta}\min\Big\{\lambda, \frac{\mu}{\sqrt{1 +\mu^2}}\Big\}
\max\{d_A(x), d_B(x)\}\\
&\geq 
\frac{\alpha\sqrt{1 -\theta}}{\kappa}\min\Big\{\lambda, \frac{\mu}{\sqrt{1 +\mu^2}}\Big\}
d_{A\cap B}(x) =\nu d_{A\cap B}(x),
\end{align}
\end{subequations}
which completes the proof.
\endproof

\section{Linear convergence of cyclic algorithms}
\label{s:lincyc}
 
We start with an elementary result.
\begin{lemma}
\label{l:linear}
Let $C$ be a closed subset of $X$, let $w\in C$, and let $(x_n)_\nnn$ be a sequence in $X$. Suppose that one of the following assumptions holds:
\begin{enumerate}
\item
\label{l:linear_seq} 
There exist $\delta \in \RPP$, $\rho\in \left[0,1\right[$, and $\sigma \in \RPP$ such that
\begin{equation}
\label{e:x+x}
\forall\nnn:\quad x_n\in \ball{w}{\delta} \ \Rightarrow\ d_C(x_{n+1}) \leq \rho d_C(x_n) \text{~and~} \|x_{n+1} -x_n\| \leq \sigma d_C(x_n).
\end{equation}
\item
\label{l:linear_op} 
The sequence $(x_n)_\nnn$ is generated by an operator $T\colon X\To X$ and there exist $\delta,\sigma \in \RPP$ and $\rho\in \left[0,1\right[$ such that
\begin{equation}
\forall x\in \ball{w}{\delta},\ \forall x_+\in Tx:\quad d_C(x_+) \leq \rho d_C(x) \text{~and~} \|x_+ -x\| \leq \sigma d_C(x).
\end{equation}
\end{enumerate}
Then if either $(x_n)_\nnn\subset \ball{w}{\delta}$ or $x_0\in \bball{w}{\frac{\delta(1-\rho)}{\sigma +1 -\rho}}$,
there exists $\ox \in C\cap \ball{w}{\delta}$ such that 
\begin{equation}\label{e:160803a}
\forall\nnn:\quad \|x_n -\ox\| \leq \frac{\sigma d_C(x_0)}{1-\rho}\rho^n,
\end{equation}
i.e., the sequence $(x_n)_\nnn$ converges $R$-linearly to a point in $C$ with rate $\rho$.
\end{lemma}
\proof 
It suffices to prove the result for \ref{l:linear_seq} because if \ref{l:linear_op} holds, then \ref{l:linear_seq} also holds for $(x_n)_\nnn$. Suppose \ref{l:linear_seq} holds, we distinguish two cases.

\emph{Case~1:} $(x_n)_\nnn\subset\ball{w}{\delta}$. Combining with \eqref{e:x+x}, we have
\begin{equation}
(\forall\nnn)\quad d_C(x_{n+1}) \leq \rho d_C(x_n) \text{~and~} \|x_{n+1} -x_n\| \leq \sigma d_C(x_n),
\end{equation} 
For each $\nnn$, take $x_n^* \in P_Cx_n$. On the one hand,
\begin{equation}
\|x_n -x_n^*\| =d_C(x_n) \leq \rho^n d(x_0) \to 0 \text{~as~} n \to +\infty.
\end{equation}
On the other hand, for all $\nnn$ and $k \in \NN\smallsetminus \{0\}$, 
\begin{equation}
\label{e:Cauchy}
\|x_n -x_{n+k}\| \leq \sum_{i=n}^{n+k-1} \|x_{i+1} -x_i\| \leq \sum_{i=n}^{n+k-1} \sigma\rho^i d(x_0) 
\leq \frac{\sigma d_C(x_0)}{1-\rho}\rho^n \to 0 \text{~as~} n \to +\infty.
\end{equation}
So $(x_n)$ is a Cauchy sequence. Therefore, $(x_n)_\nnn$ and $(x_n^*)_\nnn$ both converge to the same limit $\ox \in C\cap \ball{w}{\delta}$. 
We then obtain \eqref{e:160803a} by leting $k \to +\infty$ in \eqref{e:Cauchy}.

\emph{Case~2:} $x_0\in \bball{w}{\frac{\delta(1-\rho)}{\sigma +1 -\rho}}$. We show that this is an instance of \emph{Case~1} by proving
\begin{equation}\label{e:xnBall}
(x_n)_\nnn\subset\ball{w}{\delta}.
\end{equation}
Clearly, $\|x_0-w\|\leq \frac{\delta(1-\rho)}{\sigma +1 -\rho}\leq \delta$. So \eqref{e:xnBall} holds for $n =0$. Suppose \eqref{e:xnBall} holds for $0,1, \dots, n-1$, we shall prove that it also holds for $n$. Indeed, the induction hypothesis and \eqref{e:x+x} yield
\begin{equation}
\forall i \in \{0,1, \dots, n-1\}:\quad d_C(x_{i+1}) \leq \rho d_C(x_i) \text{~and~} \|x_{i+1} -x_i\| \leq \sigma d_C(x_i).
\end{equation} 
Noting that $d_C(x_0) \leq \|x_0 -w\|$, we obtain
\begin{subequations}
\begin{align}
\|x_n -w\| &\leq \sum_{i=0}^{n-1} \|x_{i+1} -x_i\| +\|x_0 -w\| \leq \sigma\sum_{i=0}^{n-1} d_C(x_i) +\|x_0 -w\| \\
&\leq \sigma\sum_{i=0}^{n-1} \rho^i d_C(x_0) +\|x_0 -w\| \leq \big(\sigma\frac{1}{1 -\rho} +1\big)\|x_0 -w\| \leq \delta.
\end{align}
\end{subequations}
Thus, \eqref{e:xnBall} holds for $n$. By mathematical induction principle, \eqref{e:xnBall} holds for all $\nnn$.
The conclusion now follows from \emph{Case~1}.
\endproof

\begin{corollary}
{\rm(\cite[Proposition~2.11]{Pha14})}
\label{c:linear}
Let $T\colon X\To X$ be an operator, let $C$ be a closed subset of $X$, let $w\in C$, and let $(x_n)_\nnn$ be a sequence generated by $T$. Suppose that there exist $\delta \in \RPP$ and $\rho\in \left[0, 1\right[$ such that
\begin{equation}
\forall x\in \ball{w}{\delta},\ \forall x_+\in Tx,\ \forall p\in P_Cx:\quad \|x_+ -p\| \leq \rho \|x -p\|=\rho d_C(x).
\end{equation}
Then whenever $x_0\in \ball{w}{\frac{\delta(1-\rho)}{2}}$,  
there exists $\ox \in C\cap \ball{w}{\delta}$ such that 
\begin{equation}
\forall\nnn:\quad \|x_n -\ox\| \leq \frac{(1+\rho)\|x_0-w\|}{1-\rho}\rho^n,
\end{equation}
i.e., the sequence $(x_n)_\nnn$ converges $R$-linearly to a point in $C$ with rate $\rho$.
\end{corollary}
\proof 
Let $x\in \ball{w}{\delta}$, let $x_+\in Tx$ and let $p\in P_Cx$. 
By assumption, 
\begin{equation}
d_C(x_+) \leq \|x_+ -p\| \leq \rho\|x -p\| =\rho d_C(x),
\end{equation}
and also
\begin{equation}
\|x_+ -x\| \leq \|x_+ -p\| +\|x -p\| \leq (1 +\rho)\|x -p\| =(1 +\rho)d_C(x).
\end{equation}
Now apply Lemma~\ref{l:linear}\ref{l:linear_op} with $\sigma =1 +\rho$ and note that $d_C(x_0) \leq \|x_0 -w\|$.
\endproof

The following result proves that if distance to the feasible set is reduced at least by a factor $\rho\in\left[0,1\right[$ after every {\em fixed number} of steps, then $R$-linear convergence is achieved.
\begin{lemma}[linear reduction after $k$ steps]
\label{l:redk}
Let $C$ be a closed subset of $X$, let $w\in C$, and let $(x_n)_\nnn$ be a sequence in $X$. Let also $k\in\NN\smallsetminus\{0\}$, $\delta_0\in\RPP$, $\rho\in\left[0,1\right[$, $\Gamma\in\RP$, and suppose that for every tuple $(z_0,z_1,\ldots,z_k):=(x_{kn},x_{kn+1},\ldots,x_{kn+k})$ with $z_0\in\ball{w}{\delta_0}$, we have
\begin{subequations}\label{e:redk1}
\begin{align}
&d_C(z_{k})\leq\rho d_C(z_{0})\label{e:redk1a}\quad\text{and}\\
&\forall i\in\{1,\ldots,k\},\ \forall p\in C\cap\ball{w}{2\delta_0}:\quad
\|z_{i}-p\|\leq\Gamma \|z_{0}-p\|.\label{e:redk1b}
\end{align}
\end{subequations}
Then if either $(x_{kn})_\nnn\subset\ball{w}{\delta_0}$ or $x_0\in\bball{w}{\frac{\delta_0(1-\rho)}{2+\Gamma-\rho}}$, the sequence $(x_n)_\nnn$ converges $R$-linearly to a point in $C$ with rate $\rho^{1/k}$. More specifically, there exists $\ox \in C\cap \ball{w}{\delta_0}$ such that 
\begin{equation}\label{e:160824a}
\forall\nnn:\quad \|x_n -\ox\|\leq
\frac{\Gamma(1+\Gamma)d_C(x_0)}{1-\rho}\rho^{\lfloor\frac{n}{k}\rfloor},
\end{equation}
where $\lfloor\frac{n}{k}\rfloor$ is the largest integer not exceeding $\frac{n}{k}$.
\end{lemma}
\proof
Consider the sequence $(y_n:=x_{kn})_{\nnn}$. Suppose $y_n\in\ball{w}{\delta_0}$ and take $p\in P_C x_{kn} =P_Cy_n\subset P_C(\ball{w}{\delta_0})\subseteq C\cap\ball{w}{2\delta_0}$ (see Lemma~\ref{l:ImBall}\ref{l:ImBall_P}). Then \eqref{e:redk1a} means $d_C(y_{n+1})\leq\rho d_C(y_n)$ and \eqref{e:redk1b} yields
\begin{equation}
\|y_{n+1}-y_n\|\leq
\|y_{n+1}-p\|+\|y_n-p\|
\leq \Gamma d_C(y_n) + d_C(y_n)=(1+\Gamma)d_C(y_n).
\end{equation}
So, by Lemma~\ref{l:linear}, if $(y_n)_\nnn\subset\ball{w}{\delta_0}$ or $y_0\in\bball{w}{\frac{\delta_0(1-\rho)}{2+\Gamma-\rho}}$, the sequence $(y_n)_\nnn$ converges $R$-linearly to some $\ox\in C\cap\ball{w}{\delta_0}$ and
\begin{equation}
\|y_n-\ox\|\leq\frac{(1+\Gamma)d_C(y_0)}{1-\rho}\rho^n
=\frac{(1+\Gamma)d_C(x_0)}{1-\rho}\rho^n.
\end{equation}
Now \eqref{e:redk1b} implies that for every $i\in\{1,\ldots,k\}$,
\begin{equation}
\|x_{kn+i}-\ox\|\leq\Gamma\|x_{kn}-\ox\|
=\Gamma\|y_n-\ox\|\leq\frac{\Gamma(1+\Gamma)d_C(x_0)}{1-\rho}\rho^n.
\end{equation}
Replacing $kn+i$ by $n$, we obtain
\begin{equation}
\|x_{n}-\ox\|\leq
\frac{\Gamma(1+\Gamma)d_C(x_0)}{1-\rho}\rho^{\lfloor\frac{n}{k}\rfloor}.
\end{equation}
Now if $\rho=0$, then $x_n=\ox$ for all $n\geq 1$; and if $\rho>0$, then $\rho^{\lfloor \frac{n}{k}\rfloor}\leq\frac{1}{\rho}\cdot \rho^{\frac{n}{k}}$. The lemma is proved.
\endproof

We next analyze 
the performance of $m$ steps of cyclic algorithms for quasi firmly Fej\'er monotone operators.
\begin{lemma}[consecutive steps of cyclic algorithms]
\label{l:mstep}
Let $w\in C:=\bigcap_{i\in I}C_i$, $\delta\in \RPP$, and $\nu\in \left]0, 1\right]$. For every $i\in I$, let $\gamma_i \in \left[1, +\infty\right[$ and $\beta_i\in \RPP$. Set $\O :=C\cap \ball{w}{\delta}$, $\Gamma :=(\gamma_1\cdots \gamma_m)^{1/2}$, and $\delta_0 :=\frac{\delta}{2\Gamma}\gamma_m^{1/2}$. Let $x_0, x_1, \dots, x_m$ be $m+1$ consecutive points of the cyclic algorithm with respect to $(T_i)_{i\in I}$ such that
\begin{equation}
x_0\in \ball{w}{\delta_0}
\quad\text{and}\quad
\forall i \in I:\ x_i \in T_ix_{i-1}.
\end{equation}
Then the following hold:
\begin{enumerate}
\item\label{l:mstep_qF} 
If for every $i \in I$, $T_i$ is $(\O, \gamma_i)$-quasi Fej\'er monotone on $\ball{w}{\delta/2}$, then
\begin{subequations}
\begin{align}
\label{e:xiBall}
\forall i \in I\smallsetminus \{m\}:\quad \|x_i -w\| &\leq (\gamma_1 \cdots \gamma_i)^{1/2}\delta_0\leq 
\frac{\delta}{2}, \\
\label{e:xixi-1}
\forall i \in I,\ \forall p \in \O:\quad \|x_i - p\| &\leq \gamma_i^{1/2}\|x_{i-1} -p\| \leq (\gamma_1 \cdots \gamma_i)^{1/2}\|x_0 -p\|\leq \Gamma\|x_0 -p\|.
\end{align}
\end{subequations}
\item \label{l:mstep_qfF}
If for every $i \in I$, $T_i$ is both $(\O, \gamma_i, \beta_i)$-quasi firmly Fej\'er monotone and $(C_i, \nu)$-quasi coercive on $\ball{w}{\delta/2}$, then
\begin{equation}\label{e:xmx0}
\forall p \in \O:\quad \|x_m -p\|^2 \leq (\gamma_1 \cdots \gamma_m)\|x_0 -p\|^2 -\beta\nu^2\max_{i\in I}d^2_{C_i}(x_0),\quad
\text{where $\beta :=\Big(\sum_{i\in I} \frac{1}{\beta_i}\Big)^{-1}$}.
\end{equation}
\end{enumerate}
\end{lemma}
\proof
Let $p\in \O =C\cap \ball{w}{\delta}$. 

\ref{l:mstep_qF}: First, we have $\|x_0 -w\|\leq \delta_0=\frac{\delta}{2(\gamma_1\gamma_2 \cdots \gamma_{m-1})^{1/2}} \leq \frac{\delta}{2}$ since $\gamma_i \geq 1$ for every $i\in I$. The $(\O, \gamma_1)$-quasi Fej\'er monotonicity of $T_1$ on $\ball{w}{\delta/2}$ and Lemma~\ref{l:ImBall}\ref{l:ImBall_T} then implies that
\begin{equation}
\|x_1 - p\| \leq \gamma_1^{1/2}\|x_0 -p\|
\quad\text{and}\quad
\|x_1 -w\| \leq \gamma_1^{1/2}\delta_0 = \frac{\delta}{2(\gamma_2 \cdots \gamma_{m-1})^{1/2}}
\leq \frac{\delta}{2}. 
\end{equation} 
Repeating the argument for $x_1,\ldots,x_{m-1}$, we get \eqref{e:xiBall} and the first part of \eqref{e:xixi-1}, from which the rest follows.

\ref{l:mstep_qfF}: Since quasi firm Fej\'er monotonicity implies quasi Fej\'er monotonicity, \eqref{e:xiBall} holds due to \ref{l:mstep_qF}, that is, $x_0,x_1,\ldots, x_{m-1}\in\ball{w}{\delta/2}$. Now since each $T_i$ is $(\O, \gamma_i, \beta_i)$-quasi firmly Fej\'er monotone on $\ball{w}{\delta/2}$, we derive that
\begin{subequations}
\begin{align}
\|x_1 -p\|^2 +\beta\|x_0-x_1\|^2 &\leq \gamma_1\|x_0 -p\|^2,\\
\|x_2 -p\|^2 +\beta\|x_1-x_2\|^2 &\leq \gamma_2\|x_1 -p\|^2,\\
&\vdots\\
\|x_m -p\|^2 +\beta\|x_{m-1}-x_m\|^2 &\leq \gamma_m\|x_{m-1} -p\|^2,
\end{align}
\end{subequations}
and so
\begin{subequations}
\begin{align}
\gamma_2\cdots\gamma_m\|x_1 -p\|^2 +\gamma_2\cdots\gamma_m\beta\|x_0-x_1\|^2 &\leq \gamma_1\gamma_2\cdots\gamma_m\|x_0 -p\|^2,\\
\gamma_3\cdots\gamma_m\|x_2 -p\|^2 +\gamma_3\cdots\gamma_m\beta\|x_1-x_2\|^2 &\leq \gamma_2\gamma_3\cdots\gamma_m\|x_1 -p\|^2,\\
&\vdots\\
\gamma_m\|x_{m-1} -p\|^2 +\gamma_m\beta\|x_{m-2}-x_{m-1}\|^2 &\leq \gamma_{m-1}\gamma_m\|x_{m-2} -p\|^2,\\
\|x_m -p\|^2 +\beta\|x_{m-1}-x_m\|^2 &\leq \gamma_m\|x_{m-1} -p\|^2.
\end{align}
\end{subequations}
Using the telescoping technique and the fact that $\gamma_i \geq 1$, we get
\begin{equation}
(\gamma_1 \cdots \gamma_m)\|x_0 -p\|^2
\geq \|x_m -p\|^2 +\sum_{j\in I} \beta_j\|x_{j-1} -x_j\|^2.
\end{equation}
Now the coordinate version of Cauchy--Schwarz inequality yields
\begin{equation}\label{e:x0xm}
(\gamma_1 \cdots \gamma_m)\|x_0 -p\|^2
\geq \|x_m -p\|^2 +\beta\Big(\sum_{j\in I} \|x_{j-1} -x_j\|\Big)^2.
\end{equation}
For every $i \in I$, $T_i$ is $(C_i, \nu_i)$-quasi coercive on $\ball{w}{\delta/2}$, so
$\|x_{i-1} -x_i\| \geq \nu d_{C_i}(x_{i-1})$. Hence, 
\begin{subequations}
\begin{align}
\forall i\in I:\quad
\sum_{j\in I}\|x_{j-1} -x_j\| &\geq \|x_0 -x_{i-1}\| +\|x_{i-1} -x_i\| \\
&\geq \|x_0 -x_{i-1}\| +\nu d_{C_i}(x_{i-1}) \\
&\geq \nu\big(\|x_0 -x_{i-1}\| + d_{C_i}(x_{i-1})\big)
\quad\text{(because $1\geq \nu\geq 0$)}\\
&\geq \nu d_{C_i}(x_0),
\end{align}
\end{subequations}
which yields
\begin{equation}
\sum_{j\in I}\|x_{j-1} -x_j\| \geq \nu\max_{i\in I}d_{C_i}(x_0). 
\end{equation}
Combining with \eqref{e:x0xm}, we obtain \eqref{e:xmx0}.
\endproof

The following theorems are cornerstones in our convergence analysis. In the sequel, we denote $[\rho]_+ :=\max\{0, \rho\}$ for $\rho \in \RR$.

\begin{theorem}[cyclic sequence of quasi firmly Fej\'er monotone operators]
\label{t:dist_qfF}
Let $w\in C:=\bigcap_{i\in I}C_i$, $\delta\in \RPP$, and $\nu\in \left]0, 1\right]$. For every $i\in I$, let $\gamma_i\in \left[1, +\infty\right[$ and $\beta_i\in \RPP$. Set $\O :=C\cap \ball{w}{\delta}$ and let $(x_n)_\nnn$ be a cyclic sequence generated by $(T_i)_{i\in I}$. Suppose that
\begin{enumerate}[label={\rm(\alph*)}]
\item $\{C_i\}_{i\in I}$ is $\kappa$-linearly regular on $\ball{w}{\delta/2}$ for some $\kappa \in \RPP$.
\item For every $i \in I$, $T_i$ is $(\O, \gamma_i, \beta_i)$-quasi firmly Fej\'er monotone and $(C_i, \nu)$-quasi coercive on $\ball{w}{\delta/2}$.
\end{enumerate}
Set $\Gamma :=(\gamma_1\cdots \gamma_m)^{1/2}$ and $\delta_0 :=\frac{\delta}{2\Gamma}\gamma_m^{1/2}$. 
Then 
\begin{equation}
\label{e:dxmdx0}
\forall x_0\in\ball{w}{\delta_0}:\quad d_C(x_m)\leq \rho d_C(x_0), 
\quad\text{where~} \rho :=\bigg[\Gamma^2-\frac{\nu^2}{\kappa^2}\Big(\sum_{i\in I}\frac{1}{\beta_i}\Big)^{-1}\bigg]_+^{1/2}.
\end{equation}
Consequently, if $\rho<1$ and either $(x_{mn})_\nnn\subset\ball{w}{\delta_0}$ or $x_0\in\bball{w}{\frac{\delta_0(1-\rho)}{2+\Gamma-\rho}}$, then $(x_n)_\nnn$ converges $R$-linearly to some point $\ox\in C$ with rate $\rho^{1/m}$.
\end{theorem}
\proof
Let $x_0\in \ball{w}{\delta_0}\subseteq \ball{w}{\delta/2}$. 
Since $\{C_i\}_{i\in I}$ is $\kappa$-linearly regular on $\ball{w}{\delta/2}$,
\begin{equation}
\max_{i\in I}d_{C_i}(x_0) \geq \frac{1}{\kappa}d_C(x_0).
\end{equation}
Setting $\beta :=\big(\sum_{i\in I}\frac{1}{\beta_i}\big)^{-1}$, 
Lemma~\ref{l:mstep}\ref{l:mstep_qfF} then implies that
\begin{subequations}
\label{e:xmx0'}
\begin{align}
\forall p\in \O:\quad \|x_m -p\|^2 
&\leq (\gamma_1\gamma_2 \cdots \gamma_m)\|x_0 -p\|^2 -\beta\nu^2\max_{i\in I}d^2_{C_i}(x_0) \\
&\leq \Gamma^2\|x_0 -p\|^2 -\frac{\beta\nu^2}{\kappa^2}d_C^2(x_0).
\end{align}
\end{subequations}
Letting $p\in P_Cx_0$ and noting from Lemma~\ref{l:ImBall}\ref{l:ImBall_P} that $P_Cx_0\subseteq P_C(\ball{w}{\delta/2})\subseteq \O$, we have
\begin{equation}
\label{e:4dxmdx0}
\|x_m -p\|^2 \leq \Gamma^2d_C^2(x_0) -\frac{\beta\nu^2}{\kappa^2}d_C^2(x_0) 
\leq \left[\Gamma^2 -\frac{\beta\nu^2}{\kappa^2}\right]_+ d_C^2(x_0),
\end{equation}
which leads to \eqref{e:dxmdx0}. 

Now assume $\rho<1$. Since $T_i$ is also $(\O, \gamma_i)$-quasi Fej\'er monotone on $\ball{w}{\delta/2}$ and $C\cap \ball{w}{2\delta_0}\subseteq\O$, we obtain from \eqref{e:xixi-1} in Lemma~\ref{l:mstep}\ref{l:mstep_qF} that
\begin{equation}
\forall i\in I,\ 
\forall p\in C\cap \ball{w}{2\delta_0}:\quad 
\|x_i -p\| \leq \Gamma\|x_0-p\|.
\end{equation}
By combining with \eqref{e:dxmdx0}, for every tuple $(z_0,z_1,\ldots,z_m):= (x_{mn},x_{mn+1},\ldots,x_{mn+m})$ with $z_0\in \ball{w}{\delta_0}$,  
one has  
\begin{subequations}
\begin{align}
&d_C(z_{m})\leq\rho d_C(z_{0})\quad\text{and}\\
&\forall i\in I,\ \forall p\in C\cap\ball{w}{2\delta_0}:\quad
\|z_{i}-p\|\leq \Gamma\|z_{0}-p\|,
\end{align}
\end{subequations}
which fulfills \eqref{e:redk1} with $k=m$.
The result then follows from Lemma~\ref{l:redk}.
\endproof

\begin{remark}
Regarding \eqref{e:4dxmdx0} in the proof of Theorem~\ref{t:dist_qfF}, we see that the term $\Gamma^2 -\frac{\beta\nu^2}{\kappa^2}$ is necessarily nonnegative if $x_0\notin C$; however, no general conclusion about this term can be drawn otherwise. We therefore use the notation $[\cdot]_+$
to ensure nonnegativity.
\end{remark}

Now we prove linear convergence result for cyclic sequences when there is one quasi Fej\'er monotone operator. Clearly, we need at least two operators, i.e., $m =|I|\geq 2$. Here and in what follows, $|I|$ denotes the number of elements in the set $I$.

\begin{theorem}[cyclic sequence with one quasi Fej\'er monotone operator]
\label{t:dist_qF}
Let $w\in C:=\bigcap_{i\in I}C_i$, $\delta\in \RPP$, $\nu\in \left]0, 1\right]$, and $\gamma_i\in \left[1, +\infty\right[$ for every $i\in I$. Set $\O :=C\cap \ball{w}{\delta}$ and let $(x_n)_\nnn$ be a cyclic sequence generated by $(T_i)_{i\in I}$. Suppose that
\begin{enumerate}[label={\rm(\alph*)}]
\item\label{t:dist_qF-i} $\{C_i\}_{i\in I}$ is $\kappa$-linearly regular on $\ball{w}{\delta/2}$ for some $\kappa \in \RPP$.
\item\label{t:dist_qF-ii} 
There is $j\in I$ such that for every $i \in I\smallsetminus\{j\}$, $T_i$ is $(\O, \gamma_i, \beta_i)$-quasi firmly Fej\'er monotone and $(C_i, \nu)$-quasi coercive on $\ball{w}{\delta/2}$ for some $\beta_i\in \RPP$; while $T_j$ is $(\O, \gamma_j)$-quasi Fej\'er monotone on $\ball{w}{\delta/2}$ and $T_jx\subseteq C_j$ for all $x\in\ball{w}{\delta/2}$.
\end{enumerate}
Set $\Gamma :=(\gamma_1\cdots\gamma_{m})^{1/2}$ and $\delta_0 :=\frac{\delta}{2\Gamma}\gamma_j^{1/2}$.  
Then 
\begin{equation}
\label{e:dist_qF1b}
\forall x_0\in\ball{w}{\delta_0}:\quad d_C(x_m)\leq \rho d_C(x_0),
\quad\text{where~} \rho :=\bigg[\Gamma^2 -\frac{\gamma_j\nu^2}{\kappa^2}\Big(\sum_{i\in I\smallsetminus\{j\}}\frac{1}{\beta_i}\Big)^{-1}\bigg]_+^{1/2}.
\end{equation}
Consequently, if $\rho<1$ and either $(x_{mn})_\nnn\subset\ball{w}{\delta_0}$ or $x_0\in\ball{w}{\frac{\delta_0(1-\rho)}{2+\Gamma-\rho}}$, then $(x_n)_\nnn$ converges $R$-linearly to some point $\ox\in C$ with rate $\rho^{1/m}$.
\end{theorem}
\proof
It suffices to consider only the case $j =1$ because other cases are identical up to relabeling.
Set $\delta_1 :=\frac{\delta}{2(\gamma_2 \cdots \gamma_{m-1})^{1/2}}\leq \frac{\delta}{2}$ and $\beta :=\big(\sum_{i\in I\smallsetminus\{1\}}\frac{1}{\beta_i}\big)^{-1}$.
We first claim that
\begin{equation}
\label{e:dist_qF1a}
\forall x_1 \in C_1\cap\ball{w}{\delta_1}:\quad
d_C(x_m)
\leq \left[\gamma_2 \cdots \gamma_m -\frac{\beta\nu^2}{\kappa^2}\right]_+^{1/2}d_C(x_1).
\end{equation}
On the one hand, applying Lemma~\ref{l:mstep}\ref{l:mstep_qfF} to the system $(C_i)_{i\in I\smallsetminus\{1\}}$ and $m$ consecutive points $x_1, \dots, x_m$ with $x_1\in \ball{w}{\delta_1}$, 
we deduce that
\begin{equation}
\label{e:xmx1}
\forall p \in \O:\quad \|x_m - p\|^2 
\leq (\gamma_2 \cdots \gamma_m)\|x_1 -p\|^2 -\beta\nu^2\max_{i\in I\smallsetminus \{1\}}d^2_{C_i}(x_1).
\end{equation}
On the other hand, since $d_{C_1}(x_1) =0$, the linear regularity of $\{C_i\}_{i\in I}$ yields
\begin{equation}
\label{e:linreg-x1}
\max_{i\in I\smallsetminus \{1\}}d_{C_i}(x_1)=
\max_{i\in I}d_{C_i}(x_1) \geq \frac{1}{\kappa}d_C(x_1).
\end{equation}
From \eqref{e:xmx1} and \eqref{e:linreg-x1}, letting $p \in P_Cx_1\subseteq C\cap \ball{w}{\delta}=\O$ (see~Lemma~\ref{l:ImBall}\ref{l:ImBall_P}), we obtain
\begin{subequations}
\begin{align}
d^2_C(x_m)\leq\|x_m -p\|^2 &\leq (\gamma_2 \cdots \gamma_m) \|x_1-p\|^2-\beta\nu^2\max_{i\in I}d_{C_i}^2(x_1) \\
&\leq (\gamma_2 \cdots \gamma_m) d_C^2(x_1) -\frac{\beta\nu^2}{\kappa^2} d_C^2(x_1) \\
&= \left[\gamma_2 \cdots \gamma_m -\frac{\beta\nu^2}{\kappa^2}\right]_+d_C^2(x_1),
\end{align}
\end{subequations}
which implies \eqref{e:dist_qF1a}.

Now let $x_0 \in \ball{w}{\delta_0} \subseteq \ball{w}{\delta/2}$. Then $x_1\in T_1x_0\subseteq C_1$. By applying Lemma~\ref{l:ImBall}\ref{l:ImBall_T} to $T_1$, we derive that  $x_1\in \ball{w}{\gamma_1^{1/2}\delta_0} =\ball{w}{\delta_1}$ and $d_C(x_1)\leq \gamma_1^{1/2}d_C(x_0)$. Combining these with \eqref{e:dist_qF1a}, we get \eqref{e:dist_qF1b}. The rest of the proof is exactly the same as the second part of Theorem~\ref{t:dist_qfF}.
\endproof

In the next result, we show that if the coercivity assumption is replaced by the assumption that the image of each operator $T_i$ lies in the corresponding set $C_i$, then linear reduction is obtained after $m-1$ steps (instead of $m$ steps). Thus, the rate of convergence is {\em improved}. This particular condition is satisfied for certain operators such as projectors and semi-intrepid projectors (see Section~\ref{ss:intrep}).

\begin{theorem}[refined linear convergence]
\label{t:dist_m-1}
Let $w\in C:=\bigcap_{i\in I}C_i$ and $\delta\in \RPP$. For every $i\in I$, let $\gamma_i\in \left[1, +\infty\right[$ and $\beta_i\in \RPP$. Set $\O :=C\cap \ball{w}{\delta}$ and let $(x_n)_\nnn$ be a cyclic sequence generated by $(T_i)_{i\in I}$. Suppose that
\begin{enumerate}[label={\rm(\alph*)}]
\item\label{t:dist_m-1i} $\{C_i\}_{i\in I}$ is $\kappa$-linearly regular on $\ball{w}{\delta/2}$ for some $\kappa \in \RPP$.
\item\label{t:dist_m-1ii} For every $i \in I$, $T_i$ is $(\O, \gamma_i, \beta_i)$-quasi firmly Fej\'er monotone on $\ball{w}{\delta/2}$.
\item\label{t:dist_m-1iii} For every $i\in I$ and every $x \in \ball{w}{\delta/2}$, $T_ix \subseteq C_i$.
\end{enumerate}
Set $\Gamma :=\big(\frac{\gamma_1 \cdots \gamma_m}{\min_{i\in I}\gamma_i}\big)^{1/2}$, $\delta_0 :=\frac{\delta}{2\Gamma}$, and $\rho :=\big[\Gamma^2 -\frac{1}{\kappa^2}\big((\sum_{i\in I}\frac{1}{\beta_i}) -\frac{1}{\max_{i\in I}\beta_i}\big)^{-1}\big]_+^{1/2}$.
Then
\begin{equation}\label{e:dxmdx0'}
\forall i \in I,\ \forall x_i \in C_i\cap \ball{w}{\delta_0}:\quad d_C(x_{i+m-1})\leq \rho d_C(x_i).
\end{equation}
Consequently, if $\rho<1$ and either $(x_{(m-1)n})_\nnn\subset\ball{w}{\gamma_{\max}^{-1/2}\delta_0}$ or $x_0\in\ball{w}{\gamma_{\max}^{-1/2}\cdot\frac{\delta_0(1-\rho)}{2+\Gamma-\rho}}$ where $\gamma_{\max}:=\max_{i\in I}\gamma_i$, then $(x_n)_\nnn$ converges $R$-linearly to some point $\ox\in C$ with rate $\rho^{\tfrac{1}{m-1}}$.
\end{theorem}
\proof
In addition to convention \eqref{e:cvn}, we also use
$\gamma_{mn+i} :=\gamma_i$ for $\nnn$ and $i \in I$.
For every $i\in I$, it follows from \ref{t:dist_m-1iii} that
\begin{equation}
\forall x\in\ball{w}{\delta/2},\ \forall x_+\in T_ix:\quad \|x_+ -x\|\geq d_{C_i}(x),
\end{equation}
so $T_i$ is $(C_i,1)$-quasi coercive on $\ball{w}{\delta/2}$.
Hence, all assumptions in Theorem~\ref{t:dist_qF} are fulfilled.
Now let $i\in I$ and take $m$ consecutive points $(x_i, \dots, x_{i+m-1})$ of $(x_n)_\nnn$ with $x_i\in C_i\cap \ball{w}{\delta_0}$. Then 
\begin{equation}
\|x_i -w\|\leq \delta_0 
\leq \frac{\delta}{2(\gamma_{i+1}\gamma_{i+2}\cdots \gamma_{i+m-1})^{1/2}}
\leq \frac{\delta}{2(\gamma_{i+1}\cdots\gamma_{i+m-2})^{1/2}} =:\delta_1.
\end{equation}
First, noting that $\bar\beta_i :=\big(\sum_{j\in I\smallsetminus \{i\}} \frac{1}{\beta_j}\big)^{-1}\geq \beta :=\big((\sum_{i\in I}\frac{1}{\beta_i}) -\frac{1}{\max_{i\in I}\beta_i}\big)^{-1}$ and applying claim \eqref{e:dist_qF1a} in the proof of Theorem~\ref{t:dist_qF}, we have
\begin{equation}\label{e:170803a}
d_C(x_{i+m-1})
\leq \left[\frac{\gamma_1 \cdots \gamma_m}{\gamma_i} -\frac{\bar\beta_i}{\kappa^2}\right]_+^{1/2}d_C(x_i)
\leq \left[\Gamma^2 -\frac{\beta}{\kappa^2}\right]_+^{1/2}d_C(x_i),
\end{equation}
which proves \eqref{e:dxmdx0'}.
Second, since $x_i\in\ball{w}{\delta_0}$, we derive from the quasi Fej\'er monotonicity of $T_i$'s and \eqref{e:xiBall} in Lemma~\ref{l:mstep}\ref{l:mstep_qF} that $x_{i+m-2}\in \ball{w}{\delta/2}$, which together with \ref{t:dist_m-1iii} yields 
\begin{equation}\label{e:170803b}
x_{i+m-1}\in C_{i+m-1}.
\end{equation}
Third, it follows from \eqref{e:xixi-1} in Lemma~\ref{l:mstep}\ref{l:mstep_qF} that
\begin{multline}
\label{e:xijxi}	
\hspace{2cm}\forall j\in \{1,\dots, m-1\},\ \forall p\in C\cap \ball{w}{2\delta_0}:\\
\|x_{i+j} -p\| \leq \left(\gamma_{i+1}\cdots\gamma_{i+m-1}\right)^{1/2}\|x_i -p\|\leq \Gamma\|x_i -p\|.\hspace{2cm}
\end{multline} 
Taking $p\in P_Cx_i\subseteq C\cap \ball{w}{2\delta_0}$ (due to Lemma~\ref{l:ImBall}\ref{l:ImBall_P}), we obtain
\begin{equation}\label{e:170803c}
\|x_{i+m-1} -x_i\| \leq \|x_i -p\| +\|x_{i+m-1} -p\| \leq (1 +\Gamma)\|x_i -p\| =(1 +\Gamma)d_C(x_i).
\end{equation} 
So by \eqref{e:170803a}, \eqref{e:170803b}, and \eqref{e:170803c}, we have proved that 
\begin{multline}
\label{e:extendedclaim}	
\hspace{0.5cm}\forall i \in I,\ \forall x_i \in C_i\cap \ball{w}{\delta_0}: \\
x_{i+m-1}\in C_{i+m-1},\ d_C(x_{i+m-1}) \leq \rho d_C(x_i) 
\quad\text{and}\quad \|x_{i+m-1} -x_i\| \leq (1 +\Gamma)d_C(x_i).\hspace{0.5cm}
\end{multline} 

Now assume that $\rho <1$ and that either $(x_{(m-1)n})_\nnn\subset\ball{w}{\gamma_{\max}^{-1/2}\delta_0}$ or $x_0\in\ball{w}{\gamma_{\max}^{-1/2}\cdot\frac{\delta_0(1-\rho)}{2+\Gamma-\rho}}$.
We claim that 
\begin{equation}
\label{e:subsequence}
\forall\nnn:\quad x_{(m-1)n+1}\in C_{(m-1)n+1}\cap \ball{w}{\delta_0}.
\end{equation}
Indeed, if $(x_{(m-1)n})_\nnn\subset \ball{w}{\gamma_{\max}^{-1/2}\delta_0}\subseteq \ball{w}{\delta/2}$, then \eqref{e:subsequence} holds due to \ref{t:dist_m-1iii} and Lemma~\ref{l:ImBall}\ref{l:ImBall_T}. 
If $x_0\in\ball{w}{\gamma_{\max}^{-1/2}\cdot\frac{\delta_0(1-\rho)}{2+\Gamma-\rho}}$, then by using again \ref{t:dist_m-1iii} and Lemma~\ref{l:ImBall}\ref{l:ImBall_T}, $x_1\in C_1\cap \ball{w}{\frac{\delta_0(1-\rho)}{2+\Gamma-\rho}}$, and employing \eqref{e:extendedclaim} and proceeding as in the proof of Lemma~\ref{l:linear}, we get \eqref{e:subsequence}.

Finally, we deduce from \eqref{e:xijxi}, \eqref{e:extendedclaim}, and \eqref{e:subsequence} that for every tuple $(z_0,z_1,\dots,z_{m-1}) :=(x_{(m-1)n+1}, x_{(m-1)n+2}, \dots, x_{(m-1)(n+1)+1})$,
\begin{subequations}
\begin{align}
&d_C(z_{m-1})\leq \rho d_C(z_0)\quad\text{and}\\
&\forall i\in I\smallsetminus \{m\},\ \forall p\in C\cap\ball{w}{2\delta_0}:\quad
\|z_i -p\|\leq\Gamma \|z_0 -p\|,
\end{align}
\end{subequations}
which fulfills \eqref{e:redk1} with $k=m-1$. The proof is finished by applying Lemma~\ref{l:redk}.
\endproof

\section{Applications to projection algorithms}
\label{s:apps}

\subsection{Cyclic relaxed projections}

In this section, by specializing operators $T_i$ to relaxed projectors $P_{C_i}^{\lambda_i}$, we obtain linear convergence results for the cyclic relaxed projections, one of which is possibly a {\it reflection across an injectable set}. First, we give the definition for injectability.
\begin{definition}[injectable set]
Let $C$ be a nonempty closed subset of $X$ and let $\tau \in \RP$. The set $C$ is said to be \emph{$\tau$-injectable} on a subset $U$ of $X$ if 
\begin{equation}
\forall x \in U,\ \forall p \in P_Cx:\quad
\Big[p, p +\tau\tfrac{p -x}{\|p -x\|}\Big] \subseteq C 
\end{equation}
with the convention that $\tfrac{p -x}{\|p -x\|}=0$ if $p =x$.
We say that $C$ is \emph{strictly injectable} around $w \in X$ if there exist $\tau \in \RPP$ and $\delta \in \RPP$ such that $C$ is $\tau$-injectable on $\ball{w}{\delta}$. When $C$ is $\tau$-injectable on $U =X$, we simply say that $C$ is $\tau$-injectable. When $C$ is $\tau$-injectable for all $\tau\in\RP$, we say that $C$ is $\infty$-injectable.
\end{definition}

Clearly, if $\tau>\tau'\geq 0$, then $\tau$-injectability implies $\tau'$-injectability.
To give an example of injectable sets, we recall from \cite[Section~3.2]{Gof80} that a closed convex cone $K$ of $X$ is obtuse if $-K^\ominus \subseteq K$, where $K^\ominus$ is the negative polar of $K$ defined by
\begin{equation}
K^\ominus :=\menge{y \in X}{\forall x \in K:\ \scal{x}{y} \leq 0}.
\end{equation}
The following result is a variant of \cite[Lemma~2.1(v)]{BKr04}.
\begin{proposition}
\label{p:obtuse}
Let $C$ be a translation of an obtuse cone in $X$. Then
\begin{equation}
\forall \lambda \in \left[1, +\infty\right[,\ \forall x \in X:\quad P_C^\lambda x \in C.
\end{equation} 
Consequently, $C$ is $\infty$-injectable. 
\end{proposition}
\proof
By assumption, there exist a vector $c$ and an obtuse cone $K$ in $X$ such that $C =c +K$. First, we clearly have $C+K=c+K+K=c+K=C$ since $K$ is a convex cone.

Now let $x \in X$ and set $p=P_Cx$, which is unique since $C$ is convex. It is easy to check that
\begin{equation}
p-x\in -\pnX{C}(p)\subseteq -K^\ominus\subseteq K.
\end{equation} 
So, for every $\lambda \in \left[1, +\infty\right[$, 
\begin{equation}
P_C^\lambda x=(1-\lambda)x+\lambda p= p+(\lambda-1)(p-x)\subseteq C+K=C.
\end{equation}
We therefore conclude that $C$ is $\infty$-injectable.
\endproof

We now show that injectability is a generalization of the enlargement concept, which was first defined for convex sets in \cite[Definition~2]{BIK14}.
\begin{definition}[enlargement of an arbitrary set]
\label{d:enlrg}
Given a nonempty closed subset $D$ of $X$ and $\tau \in \RP$, the \emph{$\tau$-enlargement} of $D$ is defined by the set
\begin{equation}
D_{[\tau]} :=\menge{x\in X}{d_D(x) \leq \tau} =D +\ball{0}{\tau}.
\end{equation}
\end{definition}
It is clear that $D_{[0]}=D$ and $D_{[\tau]}$ is nonempty and closed.

\begin{proposition}
\label{p:en-in}
Let $\tau \in \RP$. Then every $\tau$-enlargement is $2\tau$-injectable. In particular, every ball with radius $\tau$ is $2\tau$-injectable.
\end{proposition}
\proof
Let $C$ be a $\tau$-enlargement, say, $C =D +\ball{0}{\tau}$. 
Let $x \in X\smallsetminus C$ and let $p \in P_Cx$. There exists $q\in D$ such that $\|p-q\|\leq\tau$. It follows that
$0<d_C(x)\leq\|x-q\|\leq\|x-p\|+\|p-q\|\leq d_C(x)+\tau$.
We will show that the last two equalities happen, i.e.,
\begin{equation}\label{e:170623b}
\|x-q\|=\|x-p\|+\|p-q\|=d_C(x)+\tau.
\end{equation}
Suppose otherwise, then $\|x-q\|< d_C(x)+\tau$. Setting $z :=q+\tau\frac{x-q}{\|x-q\|}$, we have
$z\in\ball{q}{\tau}\subseteq C$ and $\|x-z\|=\|x-q\|-\tau< d_C(x)$, which is a contradiction. 
So \eqref{e:170623b} is true, which implies that $p$ lies in the segment $\left[x, q\right]$ and that $\|p -q\|=\tau$. From here, we derive that $p+2\tau\tfrac{p -x}{\|p -x\|}=p+2\tau\tfrac{q-p}{\|q-p\|}=p+2(q-p)=2q-p$ and also $\left[p-q, q-p\right]\subseteq \ball{0}{\tau}$. Hence,
\begin{equation}
\Big[p, p+2\tau\tfrac{p -x}{\|p -x\|}\Big]
=\big[p,2q-p\big]=q+\big[p-q,q-p\big]\subseteq q+\ball{0}{\tau} \subseteq C,
\end{equation}
and the conclusion follows.
\endproof

\begin{remark} 
The converse of Proposition~\ref{p:en-in} is not true. For example, consider a nontrivial obtuse cone $C$ in $\RR^2$ that is strictly contained in a halfspace. Then, for every $\tau\in\RPP$, $C$ is $\tau$-injectable but is not a $\tau$-enlargement of any subset of $\RR^2$.
\end{remark}

Enlargements emerge in several applications.
For example, the design problem in civil engineering discussed in \cite{BK15} is modeled so that all constraints are represented in the form of enlargement sets. In this case, \emph{enlargements are exactly the original constraints} of the feasibility problem. In general, one should not replace an original set by its enlargements since it may significantly change the solution of the feasibility problem. Yet there are certain cases where enlargements are actually useful. For instance, in \cite{Her75}, the image reconstruction problem is to solve a system of linear equations where constant coefficients may contain inevitable noise. Such systems may not have any exact solution. Therefore, it is reasonable to allow original equations to be only satisfied within a certain tolerance. This leads to a feasibility problem with enlargement sets. Here \emph{enlargements are replacements of the original constraints}.
In both examples, the injectability property is exploited to improve convergence.

\begin{lemma}
\label{l:inject}
Let $C$ be a nonempty closed subset of $X$ and let $w \in C$. Suppose that $C$ is strictly injectable around $w$, i.e., there exist $\tau \in \RPP$ and $\delta \in \RPP$ such that $C$ is $\tau$-injectable on $\ball{w}{\delta}$. Set $\delta' :=\min\{\tau, \delta\}$. Then
\begin{equation}
\forall \lambda \in \left[1, 2\right],\ \forall x \in \ball{w}{\delta'}:\quad P_C^\lambda x \subseteq C.   
\end{equation}   
\end{lemma}
\proof
Let $\lambda \in \left[1, 2\right]$, let $x \in \ball{w}{\delta'}\subseteq \ball{w}{\delta}$, let $x_+ \in P_C^\lambda x$, and write $x_+ =(1 -\lambda)x +\lambda p =p +(\lambda -1)(p -x)$ for some $p \in P_Cx$.
Now assume that $p \neq x$, then $0 < \|p -x\| =d_C(x) \leq \|x -w\| \leq \delta' \leq \tau$. Since $\lambda \in \left[1, 2\right]$ we have
$0 \leq \lambda -1 \leq 1 \leq \frac{\tau}{\|p -x\|}$. Combining with  
the $\tau$-injectability of $C$ on $\ball{w}{\delta}$ yields
\begin{equation}
x_+ =p +(\lambda -1)(p -x) \in \Big[p, p +\tau\frac{p -x}{\|p -x\|}\Big] \subseteq C,
\end{equation}
which finishes the proof.
\endproof

We arrive at our main results on linear convergence of cyclic relaxed projections.
\begin{theorem}[cyclic relaxed projections with at most one reflection]
\label{t:linCP}
Let $w\in C:=\bigcap_{i\in I}C_i$, $\varepsilon\in \left[0, 1\right[$, and $\delta\in \RPP$. Suppose that 
\begin{enumerate}[label={\rm(\alph*)}]
\item $\{C_i\}_{i\in I}$ is $\kappa$-linearly regular on $\ball{w}{\delta/2}$ for some $\kappa\in \RPP$.
\item $\{C_i\}_{i\in I}$ is $(\varepsilon, \delta)$-regular at $w$.
\item $\lambda_i \in \left]0, 2\right]$ for every $i \in I$ and there is at most one $\lambda_j$ equal to $2$ with the corresponding set $C_j$ being $\tau$-injectable on $\ball{w}{\delta/2}$ for some $\tau\in \RPP$. 
\item Setting $\Gamma :=\prod\limits_{i\in I} ( 1+\frac{\lambda_i\varepsilon}{1 -\varepsilon})^{1/2}$, $J :=\menge{j \in I}{\lambda_j =2}$, and $\nu :=\min\limits_{i\in I\smallsetminus J}\{1,\lambda_i\}$, 
it holds that
\begin{equation}
\label{e:linCPa}
\rho :=\left[\Gamma^2 -\frac{\nu^2}{\kappa^2}\Big(\sum_{i\in I\smallsetminus J}\frac{\lambda_i}{2-\lambda_i}\Big)^{-1}\Big(\frac{1+\varepsilon}{1 -\varepsilon}\Big)^{|J|}\right]_+^{\frac{1}{2m}} <1.
\end{equation}
\end{enumerate}
Then whenever the starting point is sufficiently close to $w$, the cyclic sequence $(x_n)_\nnn$ generated by the relaxed projections $(P_{C_i}^{\lambda_i})_{i\in I}$ converges $R$-linearly to a point $\ox \in C$ with rate $\rho$.
In particular, shrinking $\delta$ if necessary so that $\delta/2\leq \tau$, the linear convergence of $(x_n)_\nnn$ is guaranteed provided that either $(x_{mn})_\nnn\subset\ball{w}{\delta_0}$ or $x_0\in\bball{w}{\frac{\delta_0(1-\rho)}{2+\Gamma-\rho}}$, where $\delta_0 :=\frac{\delta}{2\Gamma}\min\limits_{i\in I} ( 1+\frac{\lambda_i\varepsilon}{1 -\varepsilon})^{1/2}$. 
\end{theorem}
\proof
Set $\O:=C\cap\ball{w}{\delta}$ and for every $i\in I$, set $\gamma_i :=1 +\frac{\lambda_i\varepsilon}{1 -\varepsilon}$ and $\beta_i:=\frac{2-\lambda_i}{\lambda_i}$. Then $\Gamma =(\gamma_1\cdots \gamma_m)^{1/2}$ and $\delta_0 =\frac{\delta}{2\Gamma}\min\limits_{i\in I} \gamma_i^{1/2}$.
On the one hand, for every $i\in I\smallsetminus J$, Proposition~\ref{p:rproj} implies that $P_{C_i}^{\lambda_i}$ is $\big(C_i\cap \ball{w}{\delta}, \gamma_i, \beta_i\big)$-quasi firmly Fej\'er monotone on $\ball{w}{\delta/2}$. 
On the other hand, for every $i\in I$, Proposition~\ref{p:coer-rproj} implies that $P_{C_i}^{\lambda_i}$ is $(C_i, \lambda_i)$- and therefore $(C_i,\nu)$- quasi coercive on $X$. 
We consider two cases.

\emph{Case~1:} There is no $\lambda_j$ equal to $2$, i.e, $J =\varnothing$. Noting that $\delta_0\leq \frac{\delta}{2\Gamma}\gamma_m^{1/2}$, we then apply Theorem~\ref{t:dist_qfF} to derive that 
if either $(x_{mn})_\nnn\subset\ball{w}{\delta_0}$ or $x_0\in\ball{w}{\frac{\delta_0(1-\rho)}{2+\Gamma-\rho}}$, 
the sequence $(x_n)_\nnn$ converges $R$-linearly with rate
\begin{equation}
\rho =\left[\Gamma^2
-\frac{\nu^2}{\kappa^2}\Big(\sum_{i\in I}\frac{1}{\beta_i}\Big)^{-1}\right]_+^{\frac{1}{2m}}<1.
\end{equation}

\emph{Case~2:} There is only one $\lambda_j =2$, i.e, $J =\{j\}$. 
Using Lemma~\ref{l:inject} and shrinking $\delta$ so that $\delta/2\leq \tau$, we have
\begin{equation}
\forall x \in \ball{w}{\delta/2}:\quad P_{C_j}^{\lambda_j}x \subseteq C_j.
\end{equation}
It follows from Proposition~\ref{p:rproj} that $R_{C_j} =P_{C_j}^{\lambda_j}$ is $(C_j, \gamma_j)$-quasi Fej\'er monotone on $\ball{w}{\delta/2}$ with $\gamma_j =\frac{1+\varepsilon}{1-\varepsilon}$.  
Since $\delta_0\leq \frac{\delta}{2\Gamma}\gamma_j^{1/2}$, Theorem~\ref{t:dist_qF} implies that 
if either $(x_{mn})_\nnn\subset\ball{w}{\delta_0}$ or $x_0\in\ball{w}{\frac{\delta_0(1-\rho)}{2+\Gamma-\rho}}$,  
the sequence $(x_n)_\nnn$ converges $R$-linearly with rate
\begin{equation}
\rho 
=\left[\Gamma^2 -\frac{\nu^2}{\kappa^2}\Big(\sum_{i\in I\smallsetminus J}\frac{1}{\beta_i}\Big)^{-1}\Big(\frac{1+\varepsilon}{1-\varepsilon}\Big)\right]_+^{\frac{1}{2m}}
<1.
\end{equation}
Combining the two formulas for $\rho$, we obtain \eqref{e:linCPa} and complete the proof.
\endproof

In the following, we present the convergence result with refined linear rate for cyclic over-relaxed projections. In particular, if all sets are injectable, we will obtain linear reduction after every $m-1$ steps. Therefore, the upper bound for linear rate is reduced.
\begin{theorem}[cyclic over-relaxed projections for injectable sets]
\label{t:linCP_m-1}
Let $w\in C:=\bigcap_{i\in I}C_i$, $\varepsilon\in \left[0, 1\right[$, $\delta\in \RPP$, and $\tau_i\in \RP$ for every $i\in I$. Suppose that 
\begin{enumerate}[label={\rm(\alph*)}]
\item $\{C_i\}_{i\in I}$ is $\kappa$-linearly regular on $\ball{w}{\delta/2}$ for some $\kappa\in \RPP$.
\item $\{C_i\}_{i\in I}$ is $(\varepsilon, \delta)$-regular at $w$.
\item For every $i\in I$, $\lambda_i \in \left[1, 2\right[$ and $C_i$ is $\tau_i$-injectable on $\ball{w}{\delta/2}$ with $\tau_i >0$ whenever $\lambda_i >1$.
\item Setting $\Gamma :=\max\limits_{j\in I}\prod\limits_{i\in I\smallsetminus\{j\}} (1 +\frac{\lambda_i\varepsilon}{1 -\varepsilon})^{1/2}$,  
it holds that
\begin{equation}\label{e:tlinCP_m-1a}
\rho :=\left[\Gamma^2 -\frac{1}{\kappa^2}\Big(\big(\sum_{i\in I}\frac{\lambda_i}{2-\lambda_i}\big) -\min_{i\in I}\frac{\lambda_i}{2-\lambda_i}\Big)^{-1}\right]_+^\frac{1}{2(m-1)} <1.
\end{equation}
\end{enumerate}
Then whenever the starting point is sufficiently close to $w$, the cyclic sequence $(x_n)_\nnn$ generated by the relaxed projections $(P_{C_i}^{\lambda_i})_{i\in I}$ converges $R$-linearly to a point $\ox \in C$ with rate $\rho$. 
In particular, shrinking $\delta$ if necessary so that $\delta/2\leq \min\menge{\tau_i}{\lambda_i >1}$, the linear convergence of $(x_n)_\nnn$ is guaranteed provided that either $(x_{(m-1)n})_\nnn\subset\ball{w}{\overline\delta}$ or $x_0\in\bball{w}{\frac{\overline\delta(1-\rho)}{2+\Gamma-\rho}}$, where $\overline\delta :=\frac{\delta}{2\Gamma}\big(1+\max\limits_{i\in I} \frac{\lambda_i\varepsilon}{1-\varepsilon}\big)^{-1/2}$.
\end{theorem}
\proof
We first shrink $\delta>0$ if necessary so that $\delta/2\leq \min\menge{\tau_i}{\lambda_i >1}$. For every $i\in I$, note that 
\begin{equation}\label{e:linCP_m-1a}
\forall x\in\ball{w}{\delta/2}:\quad
P_{C_i}^{\lambda_i}x\subseteq C_i.
\end{equation}
Indeed, if $\lambda_i=1$ then \eqref{e:linCP_m-1a} is automatic; and if $\lambda_i>1$, then \eqref{e:linCP_m-1a} follows from Lemma~\ref{l:inject}.

Next, define $\gamma_ i:=1+\frac{\lambda_i\varepsilon}{1-\varepsilon}$ and $\beta_i:=\frac{2-\lambda_i}{\lambda_i}$ for every $i\in I$. Then
\begin{equation}
\Gamma =\Big(\frac{\gamma_1\cdots\gamma_m}{\min_{i\in I}\gamma_i}\Big)^{1/2}
\quad\text{and}\quad
\Big(\sum_{i\in I}\frac{\lambda_i}{2-\lambda_i}\Big) -\min_{i\in I}\frac{\lambda_i}{2-\lambda_i} =\Big(\sum_{i\in I}\frac{1}{\beta_i}\Big) -\frac{1}{\max_{i\in I}\beta_i}.
\end{equation}
By Proposition~\ref{p:rproj}, for every $i\in I$, $P_{C_i}^{\lambda_i}$ is $(C_i\cap \ball{w}{\delta}, \gamma_i, \beta_i)$-quasi firmly Fej\'er monotone on $\ball{w}{\delta/2}$. Now all assumptions in Theorem~\ref{t:dist_m-1} are satisfied, hence the conclusion follows.
\endproof

\begin{remark}[refined linear rate]
\label{r:rate}
One can observe that, given the same constants $\varepsilon\in\RPP$, $\lambda_i\in\left[1,2\right[$ (which yields $\nu =1$ in Theorem~\ref{t:linCP}), and $\kappa\in\RPP$, the (upper bound) rate $\rho$ in \eqref{e:tlinCP_m-1a} is {\em smaller} than the one in \eqref{e:linCPa}. Thus, if all sets are injectable, we obtain a better upper bound for the linear rate.
\end{remark}

\begin{corollary}[refined linear convergence for cyclic projections]
\label{c:cycp}
Let $w\in C:=\bigcap_{i\in I}C_i$, $\varepsilon\in \left[0,1\right[$, and $\delta\in \RPP$. Suppose that 
\begin{enumerate}[label={\rm(\alph*)}]
\item $\{C_i\}_{i\in I}$ is $\kappa$-linearly regular on $\ball{w}{\delta/2}$ for some $\kappa\in \RPP$.
\item $\{C_i\}_{i\in I}$ is $(\varepsilon, \delta)$-regular at $w$.
\item It holds that
\begin{equation}
\rho :=\left[\frac{1}{(1 -\varepsilon)^{m-1}} -\frac{1}{(m-1)\kappa^2}\right]_+^{\frac{1}{2(m-1)}} <1.
\end{equation}
\end{enumerate}
Then whenever the starting point is sufficiently close to $w$, the cyclic sequence generated by the classical projections $(P_{C_i})_{i\in I}$ converges $R$-linearly to a point $\ox \in C$ with rate $\rho$.
\end{corollary}
\proof
Apply Theorem~\ref{t:linCP_m-1} with $\lambda_i =1$ for every $i\in I$.
\endproof

The next corollary shows that when $\{C_i\}_{i\in I}$ is a linearly regular system of superregular sets, the cyclic relaxed projections converge locally with linear rate.
\begin{corollary}[cyclic relaxed projections for superregular sets]
\label{c:linCP2}
Let $w\in C:=\bigcap_{i\in I}C_i$ and let $\lambda_i \in \left]0, 2\right]$ for every $i\in I$, where there is at most one $\lambda_i$ equal to $2$ with the corresponding $C_i$ being strictly injectable around $w$.
Suppose that the system $\{C_i\}_{i\in I}$ is linearly regular around $w$ and superregular at $w$.
Then when started at a point sufficiently close to $w$, the cyclic relaxed projection sequence generated by $(P_{C_i}^{\lambda_i})_{i\in I}$ converges $R$-linearly to a point $\ox \in C$.
\end{corollary}
\proof
Let $\varepsilon \in \left]0, 1\right[$. By assumption, there exist $\kappa \in \RPP$ and $\delta \in \RPP$ such that $\{C_i\}_{i\in I}$ is $\kappa$-linearly regular on $\ball{w}{\delta/2}$ and $(\varepsilon, \delta)$-regular at $w$. Borrowing notation from Theorem~\ref{t:linCP} and noting that $\frac{\beta\nu^2}{\kappa^2}\big(\frac{1+\varepsilon}{1-\varepsilon}\big)^{|J|}\geq \frac{\beta\nu^2}{\kappa^2} >0$ and that $1 +\frac{\lambda_i\varepsilon}{1 -\varepsilon} \to 1^+$ as $\varepsilon \to 0^+$, we choose $\varepsilon$ sufficiently small and shrink $\delta$ if necessary so that $\rho <1$. Finally, apply Theorem~\ref{t:linCP}.
\endproof

Now we turn our attention to the case of convexity in which global linear convergence is expected. 
\begin{corollary}[global linear convergence of convex cyclic relaxed projections]
\label{c:cvxCP}
Suppose that for every $i \in I$, $C_i$ is convex and that $\bigcap_{i\in I_p} C_i\cap \bigcap_{i\in I\smallsetminus I_p} \reli C_i \neq\varnothing$, where $I_p :=\menge{i \in I}{C_i \text{~is polyhedral}}$. 
Let $\lambda_i \in \left]0, 2\right]$ for every $i\in I$ and suppose that there is at most one $\lambda_i$ equal to $2$ with the corresponding $C_i$ being a translation of an obtuse cone in $X$.
Then regardless of the starting point, the cyclic relaxed projection sequence generated by $(P_{C_i}^{\lambda_i})_{i\in I}$ converges $R$-linearly to a point $\ox \in C :=\bigcap_{i\in I} C_i$. 
In particular, for every starting point $x_0 \in X$, the linear rate is
\begin{equation}
\rho :=\left[1 -\frac{\nu^2}{\kappa^2}\Big(\sum_{i\in I\smallsetminus J}\frac{\lambda_i}{2-\lambda_i}\Big)^{-1}\right]_+^\frac{1}{2m},
\end{equation}
where $J :=\menge{i \in I}{\lambda_i =2}$, $\nu :=\min\limits_{i\in I\smallsetminus J}\{1, \lambda_i\}$, and $\kappa$ is a linear regularity modulus of $\{C_i\}_{i\in I}$ on $\ball{w}{\delta/2}$ for some $\delta \in \RPP$ satisfying $\delta\geq 2d_C(x_0)$.
\end{corollary}
\proof
Let $x_0\in X$, let $\delta\in \RPP$ be such that $\delta\geq 2d_C(x_0)$, and pick $w\in C$ such that $\delta\geq 2\|x_0 -w\|\geq 2 d_C(x_0)$. Let $(x_n)_\nnn$ be the cyclic sequence generated by $(P_{C_i}^{\lambda_i})_{i\in I}$ with starting point $x_0$.
Employing \cite[Corollary~5]{BBL99}, there exists $\kappa\in \RPP$ such that $\{C_i\}_{i\in I}$ is $\kappa$-linearly regular on $\ball{w}{\delta/2}$. 
By convexity, $\{C_i\}_{i\in I}$ is $(0, \infty)$-regular at every point in $X$ (see \cite[Remark~8.2(v)]{BLPW13a}), which combined with Proposition~\ref{p:rproj} implies that for every $i\in I$, $P_{C_i}^{\lambda_i}$ is $(C_i, 1)$-quasi Fej\'er monotone on $X$. In fact, $P_{C_i}^{\lambda_i} =(1 -\frac{\lambda_i}{2})\Id +\frac{\lambda_i}{2}R_{C_i}$ is even nonexpansive due to \cite[Corollary~4.10 and Remark~4.24(i)]{BC11}. 

By Proposition~\ref{p:obtuse}, the set $C_i$ corresponding to $\lambda_i =2$, if any, is $\infty$-injectable on $X$. 
We also see that $\rho <1$ and all assumptions in Theorem~\ref{t:linCP} are therefore satisfied with $\varepsilon =0$.  
Now since $x_0\in \ball{w}{\delta/2}$, Lemma~\ref{l:ImBall}\ref{l:ImBall_T} and the $(C_i, 1)$-quasi Fej\'er monotonicity of $P_{C_i}^{\lambda_i}$ yield $(x_n)_\nnn\subset \ball{w}{\delta/2}$. Hence the proof is completed by applying Theorem~\ref{t:linCP}.
\endproof

\begin{remark}\label{r:1stRP}
When $\lambda_1 =2, \lambda_2 =\cdots =\lambda_m =1$ in Theorem~\ref{t:linCP} and Corollary~\ref{c:cvxCP}, the cyclic relaxed projections is precisely the reflection-projection algorithm, whose global convergence was studied in \cite{BKr04} with the reflection across an obtuse cone. It is worth mentioning that our results are  
the {\em first} to conclude local and global $R$-linear convergence for the reflection-projection algorithm.
\end{remark}

We finish this section by two examples showing that convergence may fail even in convex settings if there are {\it more} than one $\lambda_i$ equal to $2$ or if the strict injectability of $C_j$ corresponding to $\lambda_j =2$ is violated.

\begin{example}[failure of convergence when more than one $\lambda_i$ equal to $2$]
In $X =\RR^2$, consider two convex sets $C_1 =\RP^2$ and $C_2 =(-\RP)^2$. Then $C_1$ and $C_2$ are obtuse cones and also polyhedral sets in $X$ with $C_1\cap C_2 =\{(0, 0)\} \neq\varnothing$, hence $\{C_1,C_2\}$ is linearly regular.
It is easy to see that when started at a point $x_0 =(\zeta, \xi) \in X\smallsetminus\{(0, 0)\}$, the sequence generated by $(R_{C_1}, R_{C_2})$ does not converge since it cycles between two points $(|\zeta|, |\xi|)$ and $(-|\zeta|, -|\xi|)$.  
\end{example}

\begin{example}[failure of convergence if strict injectability is violated]
Suppose $X =\RR^2$, that $C_1 =\RR\times \{0\}$, and that $C_2 =\{0\}\times \RR$. Then $C_1$ and $C_2$ are polyhedral but not strictly injectable, and $C_1\cap C_2 =\{(0, 0)\} \neq\varnothing$. Take $x_0 =(0, \xi)$ with $\xi \in \RR\smallsetminus\{0\}$, the sequence generated by $(R_{C_1}, P_{C_2})$ cycles indefinitely between $x_0 =(0, \xi)$, $x_1 =(0, -\xi)$, $x_2 =(0, -\xi)$ and $x_3 =(0, \xi)$.
\end{example}

\subsection{Cyclic semi-intrepid projections}
\label{ss:intrep}

Cyclic intrepid projections \cite{BIK14,BK15} have found their applications in solving the feasibility problem \eqref{e:feas}, notably the road design problems \cite{BK15}.
The technique is to adjust the cyclic projections such that for every projection $P_{C_i}$, one tries to be ``more aggressive" by extrapolating into the set $C_i$ whenever possible. However, there is little incentive to ``leave" the set $C_i$, therefore, the ratio is limited to which the extrapolated point remains within the set. This idea was first used in \cite{Her75} for special polyhedra named ``strips", i.e., intersections of two halfspaces with opposite normal vectors, see also \cite{BK15,HC08}; and was later generalized in \cite{BIK14} for enlargement sets. Motivated by this, we give the definition of {\em semi-intrepid projectors}.

\begin{definition}[semi-intrepid projector to injectable sets] 
Let $\alpha \in \left[0,1\right]$, let $\tau \in \RP$, let $C$ be a $\tau$-injectable set on a given set $U$ of $X$, and let $x\in X$. The \emph{$\alpha$-intrepid projection} of $x$ into $C$ is defined by 
\begin{equation}
P^{(\alpha,\tau)}_Cx =
\left\{p +(p -x)\min\{\alpha,\tfrac{\tau}{\|p -x\|}\}\ \Big|\ p \in P_Cx\right\}
\end{equation}
with the convention that $\frac{\tau}{\|p -x\|} =0$ if $p =x$.
\end{definition}

We note that $P^{(0,\tau)}_C$ and $P^{(\alpha,0)}_C$ are just the usual projector onto $C$ and that $P^{(1,\tau)}_C$ is the original intrepid projector \cite[Definition~4]{BIK14}, see also \cite{BK15}.

\begin{proposition}
\label{p:imag-intrep}	
Let $\tau \in \RP$ and let $C$ be a $\tau$-injectable set on a given set $U$ of $X$. Then
\begin{equation}
\forall \alpha\in[0,1],\ \forall x\in U:\quad P^{(\alpha, \tau)}_C x\subseteq C.
\end{equation}
\end{proposition}
\proof
The proof is straightforward from the definition.
\endproof

\begin{proposition}[quasi firm Fej\'er monotonicity of semi-intrepid operators]
\label{p:qfF-intrep}
Let $\varepsilon\in\RP$, $\delta\in\RPP$, $\tau\in\RP$, and $\alpha\in\left[0,1\right]$. Let $C$ be a $\tau$-injectable set on $\ball{w}{\delta/2}$ and suppose that $C$ is $(\varepsilon,\delta)$-regular at $w\in C$. Then the semi-intrepid projector $P_C^{(\alpha,\tau)}$ is $(C\cap\ball{w}{\delta},\tfrac{1+\alpha\varepsilon}{1-\varepsilon}, \tfrac{1-\alpha}{1+\alpha})$-quasi firmly Fej\'er monotone on $\ball{w}{\delta/2}$.
\end{proposition}
\proof 
Take $x\in\ball{w}{\delta/2}$ and $x_+\in P_C^{(\alpha,\tau)}x$. There exists $p\in P_C x$ such that
\begin{equation}
x_+=p+\alpha'(p-x)=x+(1+\alpha')(p-x),
\quad\text{where}\quad
\alpha':=\min\{\alpha,\tfrac{\tau}{\|x-p\|}\}.
\end{equation}
Then $x_+$ is an image of the relaxed projection $P_C^{1+\alpha'}x$ and, by Proposition~\ref{p:rproj},
\begin{equation}
\forall \ox\in C\cap\ball{w}{\delta}:\quad
\|x_+-\ox\|^2+\tfrac{2-(1+\alpha')}{1+\alpha'}\|x_+-x\|^2
\leq\Big(1+\tfrac{(1+\alpha')\varepsilon}{1-\varepsilon}\Big)\|x-\ox\|^2.
\end{equation}
As $\alpha'\leq\alpha$, one can check that
$\tfrac{2-(1+\alpha')}{1+\alpha'}\geq \tfrac{1-\alpha}{1+\alpha}$ and $1+\tfrac{(1+\alpha')\varepsilon}{1-\varepsilon}\leq \tfrac{1+\alpha\varepsilon}{1-\varepsilon}$.
Hence,
\begin{equation}
\forall \ox\in C\cap\ball{w}{\delta}:\quad
\|x_+-\ox\|^2+\tfrac{1-\alpha}{1+\alpha}\|x_+-x\|^2
\leq \tfrac{1+\alpha\varepsilon}{1-\varepsilon}\|x-\ox\|^2,
\end{equation}
and the proof is complete.
\endproof

We now prove the $R$-linear convergence for the cyclic semi-intrepid projections, one of which is allowed to be the original intrepid projection \cite{BIK14,BK15}.
\begin{theorem}[cyclic semi-intrepid projections]
\label{t:cycsIP}
Let $w\in C:=\bigcap_{i\in I}C_i$, $\varepsilon \in \left]0, 1\right[$, and $\delta \in \RPP$. For every $i\in I$, let $\tau_i \in \RP$ and $\alpha_i \in \left[0, 1\right]$, where there is at most one $\alpha_j$ equal to $1$. Set $J :=\menge{j \in I}{\alpha_j =1}$ and   
\begin{equation}
\Gamma :=\Big(\frac{\gamma_1\cdots \gamma_m}{\min_{i\in I} \gamma_i^{1-|J|}}\Big)^{\frac{1}{2}},
\text{~where~} \gamma_i :=\frac{1+\alpha_i\varepsilon}{1-\varepsilon} \text{~for every~} i\in I.
\end{equation} 
Suppose that 
\begin{enumerate}[label={\rm(\alph*)}]
\item $\{C_i\}_{i\in I}$ is $\kappa$-linearly regular on $\ball{w}{\delta/2}$ for some $\kappa\in \RPP$.
\item $\{C_i\}_{i\in I}$ is $(\varepsilon, \delta)$-regular at $w$.
\item For every $i \in I$, $C_i$ is $\tau_i$-injectable on $\ball{w}{\delta/2}$.
\item It holds that
\begin{equation}
\rho :=\left[\Gamma^2 -\frac{1}{\kappa^2}\Big(\sum_{i\in I\smallsetminus J}\frac{1+\alpha_i}{1-\alpha_i} -(1-|J|)\min_{i\in I\smallsetminus J}\frac{1+\alpha_i}{1-\alpha_i}\Big)^{-1}\Big(\frac{1 +\varepsilon}{1 -\varepsilon}\Big)^{|J|}\right]_+^{\frac{1}{2(m-1+|J|)}} <1. 
\end{equation}
\end{enumerate}
Then whenever the starting point is sufficiently close to $w$, the cyclic sequence $(x_n)_\nnn$ generated by semi-intrepid projections $\big(P_{C_i}^{(\alpha_i,\tau_i)}\big)_{i\in I}$ converges $R$-linearly to a point in $C$ with rate $\rho$. In particular, the linear convergence of $(x_n)_\nnn$ is guaranteed provided that either $(x_{(m-1+|J|)n})_\nnn\subset\ball{w}{\delta'}$ or $x_0\in\bball{w}{\frac{\delta'(1-\rho)}{2+\Gamma-\rho}}$, where $\delta' :=\frac{\delta}{2\Gamma}\min\limits_{i\in I} \gamma_i^{1/2}(\max\limits_{i\in I} \gamma_i^{1/2})^{|J|-1}$.
\end{theorem}
\proof 
According to Proposition~\ref{p:imag-intrep}, for every $i\in I$,
\begin{equation}
\forall x\in \ball{w}{\delta/2}:\quad P_{C_i}^{(\alpha_i,\tau_i)}x \subseteq C_i, 
\end{equation}
and $P_{C_i}^{(\alpha_i,\tau_i)}$ is thus $(C_i, 1)$-quasi coercive on $\ball{w}{\delta/2}$. 
Next, we learn from Proposition~\ref{p:qfF-intrep} that, for $i\in I\smallsetminus J$, $P_{C_i}^{(\alpha_i,\tau_i)}$ is $(C\cap \ball{w}{\delta}, \gamma_i, \frac{1-\alpha_i}{1+\alpha_i})$-quasi firmly Fej\'er monotone on $\ball{w}{\delta/2}$ and that, for $j \in J$, $P_{C_j}^{(\alpha_j,\tau_j)}$ is $(C\cap \ball{w}{\delta}, \gamma_j)$-quasi Fej\'er monotone on $\ball{w}{\delta/2}$.

\emph{Case~1:} $J =\{j\}$. In this case, $\Gamma =(\gamma_1\cdots \gamma_m)^{1/2}$ and $\delta' =\frac{\delta}{2\Gamma}\min\limits_{i\in I} \gamma_i^{1/2}\leq
\frac{\delta}{2\Gamma}\gamma_j^{1/2}$. By Theorem~\ref{t:dist_qF},   
if either $(x_{mn})_\nnn\subset\ball{w}{\delta'}$ or $x_0\in\bball{w}{\frac{\delta'(1-\rho)}{2+\Gamma-\rho}}$, 
then $(x_n)_\nnn$ converges with $R$-linear rate
\begin{equation}
\rho =\left[\Gamma^2 -\frac{1}{\kappa^2}\Big(\sum_{i\in I\smallsetminus J}\frac{1+\alpha_i}{1-\alpha_i}\Big)^{-1}\Big(\frac{1 +\varepsilon}{1 -\varepsilon}\Big)\right]_+^{\frac{1}{2m}} <1.
\end{equation}

\emph{Case~2:} $J =\varnothing$. In this case, $\Gamma=(\frac{\gamma_1 \cdots \gamma_m}{\min_{i\in I}\gamma_i})^{1/2}$ and $\delta' =(\max\limits_{i\in I} \gamma_i)^{-1/2}\delta_0$ where $\delta_0 :=\frac{\delta}{2\Gamma}$. We get from Theorem~\ref{t:dist_m-1} that if either $(x_{(m-1)n})_\nnn\subset\ball{w}{\delta'}$ or $x_0\in\bball{w}{\frac{\delta'(1-\rho)}{2+\Gamma-\rho}}$, then $(x_n)_\nnn$ converges with $R$-linear rate 
\begin{equation}
\rho =\left[\Gamma^2 -\frac{1}{\kappa^2}\Big(\big(\sum_{i\in I}\frac{1+\alpha_i}{1-\alpha_i}\big) -\min_{i\in I}\frac{1+\alpha_i}{1-\alpha_i}\Big)^{-1}\right]_+^{\frac{1}{2(m-1)}} <1.
\end{equation}
The result follows by combining two cases.
\endproof

\begin{corollary}[global linear convergence of convex cyclic semi-intrepid projections]
\label{c:cycsIP}	
Suppose that for every $i \in I$, $C_i$ is convex and that $\bigcap_{i\in I_p} C_i\cap \bigcap_{i\in I\smallsetminus I_p} \reli C_i \neq\varnothing$, where $I_p :=\menge{i \in I}{C_i \text{~is polyhedral}}$.
Suppose also that each $C_i$ is $\tau_i$-injectable for some $\tau_i\in\RP$.
Let $\alpha_i \in \left[0, 1\right]$ for every $i\in I$ and assume there is at most one $\alpha_j$ equal to $1$. Then regardless of the starting point, the cyclic semi-intrepid projection sequence generated by $(P_{C_i}^{(\alpha_i,\tau_i)})_{i\in I}$ converges $R$-linearly to a point $\ox \in C :=\bigcap_{i\in I} C_i$. 
In particular, for every starting point $x_0 \in X$, the linear rate is
\begin{equation}\label{e:170711a}
\rho :=\left[1 -\frac{1}{\kappa^2}\Big(\big(\sum_{i\in I\smallsetminus J}\frac{1+\alpha_i}{1-\alpha_i}\big) -(1-|J|)\min_{i\in I\smallsetminus J}\frac{1+\alpha_i}{1-\alpha_i}\Big)^{-1}\right]_+^\frac{1}{2(m-1+|J|)},
\end{equation}
where $J :=\menge{i \in I}{\alpha_i=1}$  
and $\kappa$ is a linear regularity modulus of $\{C_i\}_{i\in I}$ on $\ball{w}{\delta/2}$ for some $\delta \in \RPP$ satisfying $\delta\geq 2d_C(x_0)$.
\end{corollary}
\proof
Take $x_0\in X$, $\delta\geq 2d_C(x_0)$, and choose $w\in C$ such that $\delta\geq 2\|x_0-w\|\geq 2d_C(x_0)$. Then $x_0\in \ball{w}{\delta/2}$. Let $(x_n)_\nnn$ be the cyclic sequence generated by $(P_{C_i}^{(\alpha_i,\tau_i)})_{i\in I}$ with starting point $x_0$. We observe from \cite[Corollary~5]{BBL99} that $\{C_i\}_{i\in I}$ is $\kappa$-linearly regular on $\ball{w}{\delta/2}$ for some $\kappa\in \RPP$ and from \cite[Remark~8.2(v)]{BLPW13a} that $\{C_i\}_{i\in I}$ is $(0,\infty)$-regular at every point in $X$ (due to convexity).
Note that $\rho <1$ in \eqref{e:170711a}, 
so all assumptions in Theorem~\ref{t:cycsIP} are fulfilled 
with $\varepsilon=0$.

Next, since $\{C_i\}_{i\in I}$ is $(0,\infty)$-regular at every point in $X$, Proposition~\ref{p:qfF-intrep} implies that, for every $i\in I$, $P_{C_i}^{(\alpha_i,\lambda_i)}$ is $(C_i, 1)$-quasi Fej\'er monotone on $X$. Combining with $x_0\in \ball{w}{\delta/2}$ and Lemma~\ref{l:ImBall}\ref{l:ImBall_T} gives $(x_n)_\nnn\subset\ball{w}{\delta/2}$. Now apply Theorem~\ref{t:cycsIP}.
\endproof

\subsection{Cyclic generalized Douglas--Rachford algorithm} 
In this section, we work with the index set $J :=\{1, \dots, \ell\}$, where $\ell$ is a positive integer.
For every $j \in J$, let $\lambda_j, \mu_j\in \left]0, 2\right]$, let $\alpha_j\in \left]0, 1\right[$, and let $s_j, t_j\in I$ such that $s_j\neq t_j$ and that 
\begin{equation}
\label{e:J&I}
\menge{s_j}{j \in J}\cup \menge{t_j}{j\in J} =I.
\end{equation}	
We consider the \emph{cyclic generalized Douglas--Rachford algorithm} defined by $(T_j)_{j\in J}$, where
\begin{equation}
\label{e:geDR}
\forall j \in J:\quad T_j :=(1 -\alpha_j)\Id +\alpha_j P_{C_{t_j}}^{\mu_j}P_{C_{s_j}}^{\lambda_j},
\end{equation}
and shall prove that this algorithm also possesses $R$-linear convergence properties.
It is worth noting that if each $T_j$ is a classical DR operator (i.e., $\alpha_j =1/2$, $\lambda_j =\mu_j =2$), then the cyclic generalized DR algorithm is the \emph{multiple-sets DR algorithm} \cite{BLT17}. The latter reduces to the \emph{cyclic DR algorithm} \cite{BT14} when $\ell =m$, $(s_j, t_j) =(j, j +1)$ for $j =1, \dots, m -1$, and $(s_m, t_m) =(m ,1)$; and to the \emph{cyclically anchored DR algorithm} \cite{BNP15} when $\ell =m -1$, $(s_j, t_j) =(1, j +1)$ for $j =1, \dots, m -1$.

\begin{theorem}[cyclic generalized DR algorithm]
\label{t:cDR}
Let $w\in C:=\bigcap_{i\in I}C_i$.  
Suppose that the system $\{C_i\}_{i\in I}$ is superregular at $w$ and linearly regular around $w$ and that $\{C_{s_j}, C_{t_j}\}$ is strongly regular at $w$ for every $j\in J$. Then when started at a point sufficiently close to $w$, the cyclic generalized DR sequence generated by $(T_j)_{j\in J}$ converges $R$-linearly to a point $\ox \in C$.
\end{theorem}
\proof 
Let $j\in J$ and let $\varepsilon \in \left]0, 1/3\right]$. Since $\{C_i\}_{i\in I}$ is superregular at $w$, there exists $\delta\in \RPP$ such that $C_i$ is $(\varepsilon, \sqrt{2}\delta)$-regular at $w$ for every $i \in I$. Then $C_{s_j}$ and $C_{t_j}$ are $(\varepsilon, \delta)$- and $(\varepsilon, \sqrt{2}\delta)$-regular at $w$, respectively. Using Proposition~\ref{p:qfF-DR}, $T_j$ is $(C_{s_j}\cap C_{t_j}, \gamma_j, \frac{1 -\alpha_j}{\alpha_j})$-quasi firmly Fej\'er monotone on $\ball{w}{\delta/2}$, where
\begin{equation}
\label{e:limgamma_j}
\gamma_j :=1 -\alpha_j +\alpha_j\left(1 +\frac{\lambda_j\varepsilon}{1 -\varepsilon}\right)\left(1 +\frac{\mu_j\varepsilon}{1 -\varepsilon}\right) \to 1^+ \quad\text{as}\quad \varepsilon \to 0^+.
\end{equation} 
Shrinking $\delta$ if necessary, we derive from Proposition~\ref{p:coer-DR} that $T_j$ is $(C_{s_j}\cap C_{t_j}, \nu_j)$-quasi coercive on $\ball{w}{\delta/2}$, where 
\begin{equation}
\nu_j :=\frac{\alpha_j\sqrt{1 -\theta_j}}{\kappa_j}\min\Big\{\lambda_j, \frac{\mu_j}{\sqrt{1 +\mu_j^2}}\Big\}
\quad\text{for some~} \kappa_j\in \RPP \text{~and~} \theta_j\in \left]0, 1\right[.
\end{equation}

Now by the linear regularity of $\{C_i\}_{i\in I}$, we again shrink $\delta$ if necessary and find $\kappa\in \RPP$ such that 
\begin{equation}
\forall x \in \ball{w}{\delta/2}:\quad d_C(x) \leq \kappa\max_{i\in I} d_{C_i}(x).
\end{equation} 
Since $C_{s_j}\cap C_{t_j} \subseteq C_{s_j}$ and $C_{s_j}\cap C_{t_j} \subseteq C_{t_j}$, 
\begin{equation}
\forall j \in J,\ \forall x \in X:\quad \max\{d_{C_{s_j}}(x), d_{C_{t_j}}(x)\} \leq d_{C_{s_j}\cap C_{t_j}}(x). 
\end{equation}
Noting also from \eqref{e:J&I} that
\begin{equation}
\bigcap_{j\in J}(C_{s_j}\cap C_{t_j}) =\bigcap_{i\in I}C_i =C,
\end{equation}
we conclude that the system $\{C_{s_j}\cap C_{t_j}\}_{j\in J}$ is also $\kappa$-linearly regular on $\ball{w}{\delta/2}$.

Finally, set $\nu :=\min_{j\in J} \{1, \nu_j\}$. Due to \eqref{e:limgamma_j}, we can choose $\varepsilon$ sufficiently small so that
\begin{equation}
\rho :=\left[\gamma_1 \cdots \gamma_m -\frac{\nu^2}{\kappa^2}\big(\sum_{j\in J} \frac{\alpha_j}{1 -\alpha_j}\big)^{-1}\right]_+^{1/2}<1.
\end{equation}
Thus, applying Theorem~\ref{t:dist_qfF} to $(T_j)_{j\in J}$ and the corresponding sets $(C_{s_j}\cap C_{t_j})_{j\in J}$, we obtain the $R$-linear convergence.
\endproof

We recall from Remark~\ref{r:subsystem} that the linear regularity of a system together with the strong regularity of its subsystems are less restrictive than the strong regularity of that system. This observation supports the use of our \emph{separate} assumptions on linear regularity and strong regularity in Theorem~\ref{t:cDR}.

In the case $m =2$, we obtain a generalization of \cite[Theorem~4.3]{Pha14} which proves $R$-linear convergence of the classical DR algorithm for two sets. In fact, the classical DR algorithm also converges $R$-linearly in other settings where cyclic projections may not, more details can be found in \cite{BD17,BDNP16a,BDNP16b}.
\begin{corollary}[generalized DR algorithm]
\label{c:DR2}
Let $A$ and $B$ be closed subsets of $X$ and $w\in A\cap B$. Let $\lambda, \mu \in \left]0, 2\right]$, $\alpha \in \left]0, 1\right[$, and set
\begin{equation}
T :=(1 -\alpha)\Id +\alpha P_B^\mu P_A^\lambda.
\end{equation}
Suppose that the system $\{A, B\}$ is superregular and strongly regular at $w$. 
Then when started at a point sufficiently close to $w$, the generalized DR sequence generated by $T$ converges $R$-linearly to a point $\ox\in A\cap B$. 
\end{corollary} 
\proof 
Note that strong regularity implies linear regularity (see~Fact~\ref{f:str-lin}) and apply Theorem~\ref{t:cDR} with $m =2$, $\ell =1$, and $(s_1, t_1) =(1, 2)$.
\endproof

\subsection{Affine reduction for generalized Douglas--Rachford sequences}

In this section, we extend the {\em affine reduction} scheme in \cite[Section~3]{Pha14} to generalized Douglas--Rachford sequences. Let $A$ and $B$ be nonempty closed subsets of $X$. For every $\nnn$, let $\lambda_n, \mu_n\in \left]0, 2\right]$, and $\alpha_n\in \left]0, 1\right[$. A generalized DR sequence is given by
\begin{equation}\label{e:0828a}
\forall\nnn:\quad x_{n+1}\in(1-\alpha_n)x_n+\alpha_n P_{B}^{\mu_n}P_A^{\lambda_n}x_n.
\end{equation}
We start with the following extension of \cite[Lemma~3.1]{Pha14} whose elementary proof is omitted.
\begin{lemma}
\label{l:aff}
Let $C$ be a nonempty closed subset of $X$, let $L$ be an affine subspace of $X$ containing $C$, and let $\lambda\in \RP$. Then the following hold:
\begin{enumerate}
\item 
\label{l:aff_diff}
$(\Id -P_L)P_C^\lambda =(1 -\lambda)(\Id -P_L)$.
\item
\label{l:aff_comm} 
$P_LP_C^\lambda =P_C^\lambda P_L$.
\end{enumerate}
\end{lemma}

The idea behind affine reduction for DR is to show that the \emph{shadow} of any generalized DR sequence on a certain affine subspace is again a generalized DR sequence. The next lemma provides more details.
\begin{lemma}[shadows of generalized DR sequences]
\label{l:shdwDR}
Let $L$ be an affine subspace of $X$ containing $A\cup B$ and define $y_n:=P_Lx_n$ for $\nnn$. Then the following hold:
\begin{enumerate}
\item\label{l:shdwDR-i} 
$\forall\nnn:$ $y_{n+1}\in(1-\alpha_n)y_n+\alpha_n P_{B}^{\mu_n}P_A^{\lambda_n}y_n$, i.e., $(y_n)_\nnn$ is also a generalized DR sequence.
\item\label{l:shdwDR-ii} 
$\forall\nnn:$ $x_{n+1} -y_{n+1} =\big((1-\alpha_n) +\alpha_n(1-\lambda_n)(1-\mu_n)\big)(x_n -y_n)$. 
\end{enumerate}
\end{lemma}
\proof
Let $\nnn$.

\ref{l:shdwDR-i}: Then there exist $r_n\in P_A^{\lambda_n}x_n$ and $s_n\in P_B^{\mu_n}r_n$ such that $x_{n+1} =(1-\alpha_n)x_n +\alpha_ns_n$. By Lemma~\ref{l:aff}\ref{l:aff_comm}, $P_Ls_n\in P_LP_B^{\mu_n}P_A^{\lambda_n}x_n =P_B^{\mu_n}P_A^{\lambda_n}P_Lx_n =P_B^{\mu_n}P_A^{\lambda_n}y_n$.
Since $P_L$ is an affine operator (see \cite[Corollary~3.20(ii)]{BC11}), it follows that
\begin{subequations}
\begin{align}
y_{n+1} =P_Lx_{n+1} =P_L((1-\alpha_n)x_n +\alpha_n s_n) &=(1-\alpha_n)P_Lx_n +\alpha_n P_Ls_n\\
&\in (1-\alpha_n)y_n +\alpha_n P_B^{\mu_n}P_A^{\lambda_n}y_n.
\end{align}
\end{subequations}
Hence, $(y_n)_\nnn$ is a generalized DR sequence starting at $y_0$.

\ref{l:shdwDR-ii}: Using Lemma~\ref{l:aff}\ref{l:aff_diff}, we have
\begin{subequations}\label{e:shdwDR1}
\begin{align}
s_n -P_Ls_n =(1-\mu_n)(r_n -P_Lr_n)
&=(1-\mu_n)(1-\lambda_n)(x_n -P_Lx_n)\\
&=(1-\lambda_n)(1-\mu_n)(x_n -y_n),
\end{align}
\end{subequations}
which implies that
\begin{subequations}
\begin{align}
x_{n+1} -y_{n+1} &=\big[(1-\alpha_n)x_n +\alpha_n s_n\big]
-\big[(1-\alpha_n)P_Lx_n +\alpha_n P_Ls_n\big]\\
&=(1-\alpha_n)(x_n -P_Lx_n) +\alpha_n(s_n -P_Ls_n)\\
&=\big((1-\alpha_n) +\alpha_n(1-\lambda_n)(1-\mu_n)\big)(x_n -y_n).
\end{align}
\end{subequations}
The proof is complete.
\endproof

In Corollary~\ref{c:DR2}, strong regularity of $\{A,B\}$ at $w\in A\cap B$  
is not the most general condition for $R$-linear convergence of DR sequences. Indeed, it can be relaxed to \emph{affine-hull regularity} in the sense that
\begin{equation}\label{e:affreg}
N_A(w)\cap(-N_B(w))\cap (L-w) =\{0\} \quad\text{with}\quad L :=\aff(A\cup B).
\end{equation}
This condition has been observed in \cite[Theorem~4.7]{Pha14} for the classical DR sequence ($\lambda=\mu=2$ and $\alpha=1/2$). We now continue extending such result for generalized DR sequences. For simplicity of presentation, we consider only the case of {\em constant} parameters $(\lambda_n,\mu_n,\alpha_n)\equiv(\lambda,\mu,\alpha)$.
\begin{theorem}[affine reduction for generalized DR sequences]
\label{t:0830a}
Let $A$ and $B$ be closed subsets of $X$ such that $A\cap B\neq \varnothing$, $w\in A\cap B$, and $L :=\aff(A\cup B)$. Suppose that $\{A,B\}$ is superregular and affine-hull regular at $w$. Let $(x_n)_\nnn$ be a generalized DR sequence generated by $T :=(1-\alpha)\Id +P_B^\mu P_A^\lambda$ with $\lambda, \mu\in \left]0, 2\right]$ and $\alpha\in \left]0, 1\right[$. Then the following hold:
\begin{enumerate}
\item\label{t:0830a-i} 
If $\lambda =\mu =2$, then, whenever $P_Lx_0$ is sufficiently close to $w$, the sequence $(x_n)_\nnn$ converges $R$-linearly to a point $\ox\in \Fix T$ with $P_A\ox =P_B\ox\in A\cap B$.
\item\label{t:0830a-ii} 
If either $\lambda<2$ or $\mu<2$, then, whenever $P_Lx_0$ is sufficiently close to $w$, the sequence $(x_n)_\nnn$ converges $R$-linearly to a point $\ox\in A\cap B$.
\end{enumerate}
\end{theorem}
\proof  
Define $y_n :=P_Lx_n$ for $\nnn$. By Lemma~\ref{l:shdwDR}\ref{l:shdwDR-i}, $(y_n)_\nnn\subset L$ is also a generalized DR sequence generalized by $T$. By restricting our consideration within the affine subspace $L$, affine-hull regularity \eqref{e:affreg} becomes strong regularity of $\{A,B\}$ within $L$. Thus, Corollary~\ref{c:DR2} yields that $(y_n)_\nnn$ converges $R$-linearly to a point $\oy\in A\cap B$ when $P_Lx_0 =y_0$ is sufficiently close to $w$. 

Setting $\eta :=(1-\alpha) +\alpha(1-\lambda)(1-\mu)$, we have from Lemma~\ref{l:shdwDR}\ref{l:shdwDR-ii} that 
\begin{equation}
\label{e:diff}
\forall\nnn:\quad x_n -y_n =\eta^n(x_0 -y_0).
\end{equation}

\ref{t:0830a-i}: Assume $\lambda =\mu =2$. Then $\eta =1$ and, since $(y_n)_\nnn$ converges $R$-linearly to $\oy\in A\cap B$, \eqref{e:diff} implies that $(x_n)_\nnn$ converges $R$-linearly to $\ox :=\oy +(x_0 -y_0)$. Now by \cite[Corollary~3.20(i)]{BC11}, $x_0 -y_0\in (L -L)^\perp =(L-w)^\perp$, and by \cite[Lemma~3.2]{BLPW13a}, $P_A\ox =P_A\oy =\oy$ and $P_BR_A\ox =P_B(2\oy -\ox)=P_B(\oy +y_0 -x_0) =P_B\oy =\oy$. It follows that $R_BR_A\ox =\ox$, which yields $T\ox =(1-\alpha)\ox +\alpha R_BR_A\ox =\ox$ and so $\ox\in \Fix T$.

\ref{t:0830a-ii}: Assume either $\lambda <2$ or $\mu <2$. Then $\eta <1$ and, by \eqref{e:diff}, $x_n -y_n$ converges $R$-linearly to $0$. Hence, $x_n =y_n +(x_n -y_n)$ converges $R$-linearly to $\ox =\oy\in A\cap B$.
\endproof

\begin{remark}
Theorem~\ref{t:0830a}\ref{t:0830a-ii} has never been explored before even in convex settings where one would obtain global $R$-linear convergence to the intersection; while Theorem~\ref{t:0830a}\ref{t:0830a-i} was proved in \cite[Theorem~4.7]{Pha14} for the classical DR algorithm.
With some care on the parameters, Theorem~\ref{t:0830a} can certainly be extended to the case of generalized DR iterations of the form \eqref{e:0828a} with $\disp1<\inf_{\nnn}\{\lambda_n,\mu_n\}\leq\sup_{\nnn}\{\lambda_n,\mu_n\}\leq 2$ and $\disp 0<\inf_{\nnn}\alpha_n\leq\sup_{\nnn}\alpha_n<1$.
\end{remark}

\section*{Acknowledgments}
The authors thank the referees for their valuable and constructive comments.
MND was partially supported by an NSERC (Natural Sciences and Engineering Research Council of Canada) Discovery Accelerator grant of Heinz H. Bauschke (UBC Okanagan) and by a startup research grant of the University of Newcastle. HMP was partially supported by a startup research grant of UMass Lowell. This research was partly conducted during HMP's visit at UBC Okanagan in June 2016; and HMP thanks Heinz H. Bauschke for his hospitality in Kelowna, Canada.

\end{document}